\documentclass{amsart}
\usepackage{amsfonts,amssymb}
\usepackage[cmtip,arrow]{xy}
\usepackage{pb-diagram,pb-xy}

\title{Normal differential operators and deformation theory}
\author{Paul Burchard \and Herb Clemens}

\newtheorem{thm}{Theorem}[section]
\newtheorem{lem}[thm]{Lemma}
\newtheorem{cor}[thm]{Corollary}
\newtheorem{dfn}[thm]{Definition}

\let\le=\leqslant
\let\ge=\geqslant
\let\<=\langle
\let\>=\rangle
\let\tensor=\otimes
\let\hat=\widehat
\def\bdot{\mbox{\boldmath$\cdot$}}

\def\Aa{\mathcal{A}}
\def\Dd{\mathcal{D}}
\def\Ee{\mathcal{E}}
\def\Ff{\mathcal{F}}
\def\Gg{\mathcal{G}}
\def\Hh{\mathcal{H}}
\def\Ii{\mathcal{I}}
\def\Jj{\mathcal{J}}
\def\Ll{\mathcal{L}}
\def\Oo{\mathcal{O}}
\def\Uu{\mathcal{U}}
\def\Zz{\mathcal{Z}}

\DeclareSymbolFont{euscriptbold}{U}{eus}{b}{n}
\DeclareSymbolFontAlphabet\mathcursivebold{euscriptbold}
\def\HHh{\mathcursivebold{H}}

\def\CC{\mathbb{C}}

\def\HHH{\mathbf{H}}

\def\Hom{\mathop{\mathrm{Hom}}}
\def\sheafHom{\mathop{\mathfrak{Hom}}}
\def\sgn{\mathop{\mathrm{sgn}}}

\def\wideD{\setlength{\dgARROWLENGTH}{2.5em}}
\def\thinD{\setlength{\dgARROWLENGTH}{1.5em}}

\begin{document}

\begin{abstract}
This paper develops the theory of a sheaf of normal differential
operators to a submanifold $Y$ of a complex manifold $X$ as a generalization of
the normal bundle. We show that the global sections of this sheaf play an
analogous role for formal deformations of $Y$ to the role played by the normal
bundle with respect to first-order deformations.
\end{abstract}

\maketitle

\section{Introduction and Motivation}\label{sec.intro}

Our aim in this paper is to develop a more concrete and geometrically intuitive formulation of the higher-order deformation theory of subvarieties.
Our basic strategy is to use local-formal duality to reformulate the standard theory in terms of {\em local cohomology sheaves}.
These sheaves are rich in geometric meaning,
as demonstrated by their alternative interpretations as sheaves of {\em residue operators\/} or {\em normal differential operators}.
Because these sheaves already carry higher-order information,
it is possible to work more concretely (doing more at the level of sheaves, or $H^0$).

The geometric setting in which our local cohomology sheaves become functorial is that of {\em deformation categories\/} of subvarieties.
In the simplest such category, an object is a smooth subvariety $Y \arrow{e,V} X$,
and a morphism from $T \arrow{e,V} S$ to $Y \arrow{e,V} X$ is a deformation diagram
\begin{equation*}
\wideD
\begin{diagram}
   \node{Y} \arrow{s} \arrow{e,V} \node{W} \arrow{s} \arrow{e} \node{X} \\
   \node{T} \arrow{e,V} \node{S}
\end{diagram}
\end{equation*}
where the square is a pullback and the maps at $S$ are transversal.
Composition of morphisms is pullback.

Our main technical result---a special case of which is proved in~\S\ref{sec.defcrit} of this paper---is
an {\em equivalence of categories\/}
between a certain deformation category of smooth formal subvarieties,
and a corresponding category of local cohomology sheaves.
In order to obtain such an equivalence, however,
the natural homotopy product structure which these special local cohomology sheaves carry must be respected.

Let us describe the local cohomology sheaves a bit more explicitly,
if only to indicate why our global-to-formal-to-local strategy doesn't destroy the cohomological information required for deformation theory.
This may not be obvious because, for a smooth subvariety $Y \arrow{e,V} X$ of dimension $d$ and codimension $c$,
the interesting cohomology group $H^{2d+1}(\CC_X)$ is immediately killed upon restriction to the formal neighborhood $\hat{X}_Y$.
Fortunately, the Abel-Jacobi map involves only a specific Hodge level of this cohomology,
and once we restrict our attention to this level,
the relevant classes manage to sneak safely across into local cohomology:
\begin{equation*}
\thinD
\begin{diagram}
  \node{ {H^{2d+1} (\CC_{X})}^* } \arrow{e,A}
  \node{ {\HHH^{d} (\Omega_{X}^{> d})}^* }
  \node{ {\HHH^{d} (\Omega_{\hat{X}_Y}^{> d})}^* } \arrow{w} \arrow{e,t}{\delta^*}
  \node{ {\HHH^{2d} (\Omega_{\hat{X}_Y}^{\le d})}^* } \arrow{e,=}
  \node{ \HHH_{Y}^c (\Omega_X^{\ge c}) }
\end{diagram}
\end{equation*}
Thus, the ``local Hodge sheaf'' $\HHh_{Y}^c (\Omega_X^{\ge c})$
is sufficient to express Abel-Jacobi data.

Our ultimate object of study is the local cohomology of the whole ambient de~Rham complex
\begin{equation}
\label{eq.dperp}
	\left( \strut \Hh_{Y}^c (\Omega_X^{\bdot}) ,\ \partial \right).
\end{equation}
To properly handle relative subvarieties $Y/Y' \arrow{e,V} X/X'$,
we must further endow this complex with a filtration by base degree of closed forms:
\begin{equation}
\label{eq.dperprel}
	F^{p/p'} \Hh_{Y}^c (\Omega_X^{p})
	:= \HHh_{Y}^c \left( (\Omega_X^{\cdot} \wedge \Omega_{X'}^{\ge p'})^{\ge p} \right) .
\end{equation}
For the purposes of this paper, however,
we only need a particular piece of this complex---one which
can be defined directly in terms of relative forms~(\S\ref{sec.ndo}).

The special local cohomology sheaf~(\ref{eq.dperp}) carries a natural homotopy product,
which can be formalized using the notion of {\em co-implicial\/} objects.
Recall that co-simplicial objects are representations of the category of finite ordinals and ordered maps;
{\em co-implicial\/} objects are representations of the richer category which includes partially-defined ordered maps.
A co-implicial object will therefore carry {\em co-barycenter\/} maps in addition to the usual {\em co-face\/} and {\em co-degeneracy\/} maps.
These maps can be thought of as operations on product expressions implementing (respectively) factor deletions, unit insertions, and pairwise factor associations.
To avoid extra formalism, however, this paper will take a more ad-hoc approach to the product.

Our main applications will come from studying the pairings between
ambient global cohomology and the cohomology of normal differential operators.
These pairings essentially come from the local-formal pairings which were used to construct the theory in the first place;
for instance,
\begin{equation*}
\thinD
\begin{diagram}
  \node{ H^q (\HHh_{Y}^c (\Omega_X^{\ge c})) \tensor \HHH^r (\Omega_X^{\le d})  } \arrow{e}
  \node{ \HHH_{Y}^{q+c} (\Omega_X^{\ge c}) \tensor \HHH^r (\Omega_{\hat{X}_Y}^{\le d}) } \arrow{e,t}{\text{local-formal}}
  \node{ H^{q+r} (\CC_Y) }
\end{diagram}
\end{equation*}
By looking at similar pairings in the relative setting,
we are able to generalize a result of Bloch-Ran-Kawamata on the {\em semi-regularity map}.
Our generalization, given at the end of \S\ref{sec.obst},
produces an {\em obstruction differential\/} which is exact if there is no cohomological obstruction to deformation.

The authors wish to thank Henryk Hecht, J\'{a}nos Koll\'{a}r, and Dragan
Milicic of the University of Utah, Morihiko Saito of RIMS, Japan, and Claire
Voisin of the CNRS, France, for valuable conversations and extremely useful
guidance during the preparation of this paper.

\medskip
Paul Burchard $\<$burchard@pobox.com$\>$

Herbert Clemens $\<$clemens@math.utah.edu$\>$

June, 1998

\newpage

\section{Relation to the Kuranishi Theory}\label{sec.kuranishi}

To explain the results and methods of this paper in more detail, we must
first recall, at least in rough form, some basic facts from the Kuranishi
theory of deformations of complex manifold as developed by Deligne and
Goldman-Millson in~\cite{bib.GM}. Given a smooth family
\begin{equation*}
Z/Z'=\left\{ Z_{z'}\right\} _{z'\in
Z'}
\end{equation*}
of connected complex manifolds over a polydisc $Z$ with coordinate
\begin{equation*}
z'=\left( z_{1}',\ldots ,z_{n}'\right)
\end{equation*}
and distinguished point $0$, one can associate a transversely holomorphic
trivialization, that is, a $C^{\infty }$-trivialization
\begin{equation}
\begin{diagram}
\node{Z} \arrow{s,r}{ p_{Z'}} \arrow{e,t}{\left( \sigma ,p_{Z'}\right) } \node{Z_{0}\times Z'} \arrow{s} \\
\node{Z'} \arrow{e,=} \node{Z'} 
\end{diagram}
\label{eq.01}
\end{equation}
for which the fibers are all complex holomorphic polydisks. (See, for
example,~\cite{bib.C1}.) Following~\cite{bib.GM}, Kuranishi theory then associates a formal
deformation as follows (see also~\cite{bib.C}).\bigskip

\begin{lem}\label{lem.02}
The transversely holomorphic trivialization $\sigma $
distinguishes a subset of the complex-valued $C^{\infty }$-functions on
$Z$, namely those functions $f$ which restrict to a holomorphic function on
each fiber of $\sigma $ (i.e. functions with a power-series representation
\end{lem}

\begin{equation}
\sum\nolimits_{I}f_{I}z^{\prime I}:=\sum\nolimits_{I\in
\left( \Bbb{Z}^{\ge 0}\right) ^{n'}}\left(
f_{I}\circ \sigma \right) z^{\prime I}
\label{eq.021}
\end{equation}
where the $f_{I}$ are $C^{\infty }$-functions on $Z_{0}$).
Furthermore there exist global $C^{\infty }$-bundle-sections
\begin{equation*}
\xi _{J}\in \Aa^{0,1}\left( T_{Z_{0}}\right) ,\quad J\in
\left( \left( \Bbb{Z}^{\ge 0}\right) ^{n'}-\left\{
\left( 0,\ldots ,0\right) \right\} \right)
\end{equation*}
such that a function~(\ref{eq.021}) is holomorphic if and only if
\begin{equation*}
\bar{D}_{\sigma }\left( f\right) :=\left( \bar{\partial }-\sum\nolimits_{J}z^{\prime J}\xi _{J}\right) \left( f\right) =
\sum\nolimits_{I}\bar{\partial }f_{I}z^{\prime I}-\sum\nolimits_{I,J}\xi _{J}\left( f_{I}\right) z^{\prime I+J}=0
.
\end{equation*}
If we let
\begin{equation*}
\tau _{i}
\end{equation*}
denote the unique lifting of $\frac{\partial }{\partial z_{i}'}$ to
a $\left( 1,0\right) $-vector-field on $Z$ determined by the trivialization~(\ref{eq.01})
we have that as operators on functions on $Z_{0}$, for example
\begin{equation*}
\xi _{K}=-\left( \left. {}\right| _{Z_{0}}\right) \circ
\left[ L_{\tau }^{K},\bar{\partial }_{Z}\right]
\end{equation*}
for $K=\left( 0,\ldots ,0,1,0,\ldots ,0\right) .$ For
\begin{equation*}
\xi :=\sum\nolimits_{J}z^{\prime J}\xi _{J}
\end{equation*}
the integrability condition
\begin{equation}
\bar{\partial }\xi -\frac{\left[ \xi ,\xi \right] }{2}=0
\label{eq.024}
\end{equation}
is satisfied. If two transversely holomorphic trivializations of $%
Z/Z'$ differ by a holomorphic automorphism defined over $Z' $ and equal to the identity on $Z_{0}$, they give the same Kuranishi
datum $\xi $.

Continuing to follow~\cite{bib.GM}, we define the differential graded Lie algebra
\begin{equation*}
\left( \sum\nolimits_{i\ge 0}\Ll^{i},\bar{\partial }\right)
\end{equation*}
where
\begin{equation*}
\Ll^{i}=\left\{ \lambda =\sum\nolimits_{J}z^{\prime J}\lambda
_{J}:\lambda _{J}\in \Aa_{Z_{0}}^{0,i}\left(
T_{Z_{0}}\right) \right\} .
\end{equation*}
The Lie algebra of vector fields $\Ll^{0}$ acts via an exponential map on
the subspace $\breve{\Ll}^{1}$ of $\Ll^{1}$ consisting of those
elements which satisfy the integrability condition~(\ref{eq.024}). The quotient is
the space of formal deformations of $Z_{0}$ parametrized by $Z'$.

Next suppose that we have a commutative diagram
\begin{equation}
\begin{diagram}
\node{Y} \arrow{s,r}{ p_{Y'}} \arrow{e,t}{h} \node{Z} \arrow{s,r}{ p_{Z'}} \\
\node{Y'} \arrow{e,t}{h'} \node{Z'} 
\end{diagram}
\label{eq.031}
\end{equation}
of complex manifolds where $p_{Y'}$ is also a smooth
surjection and $h'$ factors as a smooth surjection followed by a closed immersion
\begin{equation*}
\begin{diagram}
\node{Y'} \arrow{e,t}{p_{Y^{\prime \prime }}} \node{Y^{\prime \prime }} \arrow{e,t}{i^{\prime \prime }} \node{Z'}
\end{diagram}
\end{equation*}
so that we have fibered product:
\begin{equation*}
\begin{diagram}
\node{Y'} \arrow{s,r}{ p_{Y^{\prime \prime }}} \arrow[2]{e,t}{i':=\left( y',h'\left( y'\right) \right) } \arrow{e,!} \node{\strut} \arrow{e,!} \node{X':=Y^{\prime }\times Z'}
	\arrow{s,r}{ p_{X^{\prime \prime }}:=\left( p_{Y^{\prime \prime }}\left( y'\right) ,z'\right) } \\
\node{Y^{\prime \prime }} \arrow[2]{e,t}{i^{\prime \prime }:=\left( y^{\prime \prime },i^{\prime \prime }\left( y^{\prime \prime }\right) \right) } \node[2]{X^{\prime \prime }:=Y^{\prime \prime }\times Z^{\prime }} 
\end{diagram}
\end{equation*}
and
\begin{equation*}
\begin{diagram}
\node{X:=X'\times _{Z'}Z} \arrow{s,r}{ p_{X'}} \arrow{e,t}{j} \node{Z} \arrow{s,r}{ p_{Z'}} \\
\node{X'=Y'\times Z'} \arrow{e,t}{j'} \node{Z'} 
\end{diagram}
\end{equation*}

Suppose that, in the induced diagram
\begin{equation}
\begin{diagram}
\node{Y} \arrow{s,r}{ p_{Y'}} \arrow{e,t}{i} \node{X} \arrow{s,r}{ p_{X'}} \arrow{e,t}{j} \node{Z} \arrow{s,r}{ p_{Z'}} \\
\node{Y'} \arrow{e,t}{i'} \node{X'} \arrow{e,t}{j'} \node{Z'} 
\end{diagram}
,
\label{eq.041}
\end{equation}
$i$ is a closed imbedding and $p_{Y'}$ is proper. For many
algebro-geometric purposes, it is desireable to characterize a certain class
of deformation/extensions of~(\ref{eq.041}), namely to characterize commutative
diagrams
\begin{equation}
\begin{diagram}
\node{Y} \arrow{s,r}{ p_{Y'}} \arrow{e}  \node{W} \arrow{s,r}{ q_{X'}} \arrow{e}  \node{X} \arrow{s,r}{ p_{X'}} \\
\node{Y'} \arrow{e,t}{i'} \node{X'} \arrow{e,=} \node{X'} 
\end{diagram}
\label{eq.042}
\end{equation}
with $q_{X'}$ a smooth surjection such that the left-hand
square is a fibered product. (If~(\ref{eq.031}) itself is a fibered product and vertical maps are
proper, there is a unique such (proper) extension.)  From~(\ref{eq.042}) we can
construct a smoothly varying family of transversely holomorphic
trivializations
\begin{equation*}
\begin{diagram}
\node{\left\{ y'\right\} \times Z} \arrow{e,t}{\sigma _{y'}} \node{Z_{h'\left( y'\right) }}
\end{diagram}
\end{equation*}
for $X/X'$ as above but with the added property that
\begin{equation*}
\sigma _{y'}^{-1}\left( Y_{y'}\right) =q_{X'}^{-1}\left( \left\{ y'\right\} \times Z'\right) .
\end{equation*}
This translates to the condition
\begin{equation*}
\xi _{y',J}\in \ker \left( \Aa_{Z_{h'\left( y'\right) }}^{0,1}\left( T_{Z_{h'\left( y'\right)
}}\right) \arrow{e} \Aa_{Y_{y'}}^{0,1}\left(
N_{Y_{y'}\backslash Z_{h'\left( y'\right)
}}\right) \right)
\end{equation*}
on the Kuranishi data in Lemma~\ref{lem.02}, where $N_{Y_{y'}\backslash Z_{h'\left( y'\right) }}$
denotes the normal bundle.

For each $y'\in Y'$, we define the differential graded Lie
algebra
\begin{equation*}
\left( \sum\nolimits_{i\ge 0}\Ll_{y'}^{i},\bar{\partial }\right)
\end{equation*}
where
\begin{equation*}
\Ll_{y'}^{i}=\left\{ \lambda =\sum\nolimits_{J}z'\lambda _{J}\left( y'\right) :\lambda _{J}\left(
y'\right) \in \ker \left( \Aa_{Z_{h'\left( y'\right) }}^{0,i}\left( T_{Z_{h'\left( y'\right)
}}\right) \arrow{e} \Aa_{Y_{y'}}^{0,i}\left(
N_{Y_{y'}\backslash Z_{h'\left( y'\right)
}}\right) \right) \right\} .
\end{equation*}
For a fixed transversely holomorphic trivialization of the deformation $%
Z/Z'$ given in~(\ref{eq.01}), we have Kuranishi data
\begin{equation*}
\xi _{z'}\in \breve{\Ll}_{z'}^{1}\subseteq \Ll_{z'}^{1}
\end{equation*}
for each $z'\in Z'$. As above the vector fields $%
\Ll_{z'}^{0}$ act on $\Ll_{z'}^{1}$ as does the subspace $%
\Ll_{y'}^{0}$ whenever $h'\left( y'\right)
=z'$. If $\xi _{h'\left( y'\right) }\in
\Ll_{y'}^{1}$ and $\alpha \in \Ll_{y'}^{0}$, then the
infinitesimal action of $\alpha $ on $\xi $ is given by the formula
\begin{equation*}
\dot{\xi}=\bar{\partial }\alpha +\left[ \alpha ,\xi \right] ,
\end{equation*}
so, if $\alpha \in \Ll_{y'}^{0}$, then $\dot{\xi}\in \Ll_{y'}^{1}$. The orbits of the exponential action of $\Ll_{y'}^{0}$ on $%
\breve{\Ll}_{y'}^{1}$ form the space of formal deformations of $%
Y_{y'}\subseteq Z_{h'\left( y'\right) }$
inside the fixed deformation $Z/Z'$.

The Kuranishi datum is invariant under holomorphic automorphisms over $%
Z'$ which restrict to the identity at time zero (see the final
assertion of Lemma~\ref{lem.02}). So, if $Y$ in~(\ref{eq.01}) is a small neighborhood of a
point, the action of $\Ll_{y'}^{0}$ on $\breve{\Ll}_{y'}^{1}$
is transitive. Thus, in the above formulation, we do not see the local
structure of the space of deformations of a small neighborhood of a point in
$Y$. Our purpose is to provide a formalism for that local structure. To
first order, that structure is classical. For example, if $Y'=\left\{ 0\right\} $, the local structure is as follows. The tangent sheaf
$T_{Z}$ of the total space $Z$ has distinguished subsheaves, which we
denote as $T_{Z/Z'},$ consisting of those derivations which
annihilate $p_{Z'}^{-1}\left( \Oo_{Z'}\right) $ where $%
p_{Z'}:Z\rightarrow Z'$ is the natural projection
associated with the family, and $\hat{T}_{Z/Z'},$ consisting
of those derivations which are ``lifts'' of derivations on $Z'$. The latter
forms a $p_{Z'}^{-1}\left( \Oo_{Z'}\right) $-module
and contains the subsheaf $T_{Z/Z'}$ of the $p_{Z'}^{-1}\left( \Oo_{Z'}\right) $-linear derivations.
($T_{Z/Z'}$, on the other hand, is an $\Oo_{Z}$-submodule of $%
T_{Z}$.) If $Y_{0}$ is a compact, connected complex $d$-dimensional
submanifold of a fiber $Z_{0}$ of the family, then all
global sections of the normal bundle $N_{Y_{0}\backslash Z}$ actually
lie in a sub-sheaf
\begin{equation*}
\hat{N}_{Y_{0}\backslash Z}=\frac{\left. \hat{T}_{Z/Z'}\right| _{Y_{0}}}{T_{Y_{0}}}.
\end{equation*}
Morphisms
\begin{equation}
T_{Z',0}\rightarrow H^{0}\left( \hat{N}_{Y_{0}\backslash
Z}\right)
\label{eq.05}
\end{equation}
give the first-order deformations of $Y_{0}$ in $Z$ with parameter
space $Z'$, and $\hat{N}_{Y_{0}\backslash Z}$ sits in an exact
sequence
\begin{equation*}
0\arrow{e} N_{Y_{0}\backslash Z_{0}}\arrow{e} \hat{N}%
_{Y_{0}\backslash Z}\arrow{e} \left. p_{Z'}^{-1}T_{Z'}\right| _{Y_{0}}\arrow{e} 0.
\end{equation*}

In this paper, we will develop analogous machinery for higher-order
deformations of $Y/Y'$ in $Z/Z'$. This analogue comes from
the theory of $\Dd$-modules, and is similar to the sheaf of ``normal
differential operators,'' the push-forward of the structure sheaf of $%
Y_{y'}$ in the sense of $\Dd$-modules. Although the sheaf
of Lie algebras
\begin{equation*}
\Ll_{Z}
\end{equation*}
given by vector fields on Z does not have a natural action on the relative
deRham complex $\Omega _{Z/Z'}^{*}$, the sub-Lie-algebra
\begin{equation*}
\Ll_{Z/Z'}
\end{equation*}
given by sections of $\hat{T}_{Z/Z'}$ does, via Lie
differentiation. Under this action, sections of $\hat{N}_{Y_{0}\backslash
Z_{0}}$ do not give well-defined operators at points of $Y_{y'}$, but rather should be thought of as equivalence classes of
operators which give well defined morphisms
\begin{equation*}
R^{j}\left( p_{Z'}\right) _{*}\left( \Omega
_{Z/Z'}^{*}\right) \arrow{e,t}{}H^{j}\left(
\Omega _{Y/Y^{^{\prime \prime }}}^{*}\right) .
\end{equation*}

These considerations lead us to use $\Dd$-modules constructed from
differential operators on the relative deRham complex of $Z/Z'$
with values in the relative deRham complex of $Y/Y^{\prime \prime }$. The
fundamental object of study is a certain sheaf $\Dd_{Y/Y^{\prime \prime
}}^{Z/Z'}$ supported on $Y$. The fiber of its pushforward $\left.
\left( \left( p_{Y'}\right) _{*}D_{Y/Y^{\prime \prime
}}^{Z/Z'}\right) \right| _{y'}$ should be thought of
as the higher-order analogue of the space
\begin{equation*}
\left( p_{Y'}\right) _{*}\hat{N}_{Y\backslash Y'\times Z}
\end{equation*}
where $\hat{N}_{Y\backslash Y'\times Z}$ sits in the exact
sequence
\begin{equation*}
0\arrow{e} N_{Y\backslash Y'\times _{Z'}Z}\arrow{e} \hat{N}_{Y\backslash Y'\times
Z}\arrow{e} p_{Y'}^{-1}{h'}^{*}T_{Z'}\arrow{e} 0.
\end{equation*}
The general construction and its basic properties and filtrations occupy
\S\S\ref{sec.dmods}--\ref{sec.functor}. The pairing of global sections with ambient cohomology is given in \S\ref{sec.resol}.

To illustrate, we return to the case in which $Y'=\left\{ 0\right\} $ is a
point. ``Almost-multiplicative'' morphisms from the differential operators
at $0\in Z'$ with values in $\Bbb{C}$ to global sections of
$\Dd_{Y_{0}/\left\{ 0\right\} }^{Z/Z'}$ will characterize
deformations of $Y_{0}$ in $Z$ with parameter space $Z'$ just
as morphisms~(\ref{eq.05}) characterize first-order deformations. This equivalence,
in its more general, relativized version, is developed in \S\S\ref{sec.gauss}--\ref{sec.defcrit}.

As indicated in \S\ref{sec.intro}, a value of this type of approach to higher-order deformation theory comes
from the fact that it is well suited for computing the interaction of
higher-order obstructions with the cohomology of the ambient manifold $%
Z/Z'$. To motivate this assertion, suppose that we are given a
smooth family $Z/Z'$ of deformations of $Z_{0}$, and let
\begin{equation*}
\begin{diagram}
\node{T_{Z',0}} \arrow{e,t}{\upsilon} \node{H^{0}\left( N_{Y_{0}\backslash Z_{0}}\right)}
\end{diagram}
\end{equation*}
be the obstruction map to deforming $Y_{0}$ to first order with $Z_{0.}$.
Recall the easy first-order fact that the following diagram is commutative:
\begin{equation*}
\setlength{\dgARROWLENGTH}{1.5em}
\begin{diagram}
\node{\wedge ^{r}H^{0}\left( N_{Y_{0}\backslash Z_{0}}\right) \tensor T_{Z',0}\tensor H^{d-1}\left( \Omega _{Z_{0}}^{d+r+1}\right)} \arrow{s,r}{\text{Kodaira-Spencer}} \arrow{e,t}{\upsilon}
	\node{H^{0}\left( \wedge ^{r}N_{Y_{0}\backslash Z_{0}}\right) \tensor H^{0}\left( N_{Y_{0}\backslash Z_{0}}\right) \tensor H^{d-1}\left( \Omega_{Z_{0}}^{d+r+1}\right)}
	\arrow{s,r}{\text{contraction}} \\
\node{\wedge^{r}H^{0}\left( N_{Y_{0}\backslash Z_{0}}\right) \tensor H^{d}\left(\Omega _{Z_{0}}^{d+r}\right)}
	\arrow{e} \node{H^{d}\left( \omega_{Y_{0}}\right)}
\end{diagram}
\end{equation*}
The contraction map across the bottom of this diagram is the infinitesimal
Abel-Jacobi mapping. The diagram tells us that, at least when the natural
mapping
\begin{equation*}
\begin{diagram}
\node{\wedge ^{r}H^{0}\left( N_{Y_{0}\backslash Z_{0}}\right) \tensor H^{d-1}\left( \Omega _{Z_{0}}^{d+r+1}\right)} \arrow{s} \\
\node{H^{0}\left( \wedge ^{r}N_{Y_{0}\backslash Z_{0}}\right) \tensor H^{d-1}\left( \wedge ^{r+1}N_{Y_{0}\backslash Z_{0}}^{*}\tensor \omega _{Y_{0}}\right)} \arrow{s} \arrow{e,!} \node{\hspace{2em}} \\
\node{H^{d-1}\left( N_{Y_{0}\backslash Z_{0}}^{*}\tensor \omega _{Y_{0}}\right)} \arrow{e,t,=}{\text{Serre duality}} \node{H^{1}\left( N_{Y_{0}\backslash Z_{0}}\right)^{*}}
\end{diagram}
\end{equation*}
is surjective, the obstruction to moving $Y_{0}$ with $Z_{0}$ is
completely computed from the infinitesimal Abel-Jacobi mapping for $%
Y_{0}$ in $Z_{0}$. A striking example is the case
\begin{equation*}
d=r=1\quad n=3\quad \omega _{Z_{0}}=\Oo_{Z_{0}},
\end{equation*}
where the surjectivity is automatic. In this case, the above diagram tells
us that the curve $Y_{0}$ deforms to first order in a direction in $%
Z'$ if and only if the corresponding element of
\begin{equation*}
H^{1}\left( \Omega _{Z_{0}}^{2}\right) =H^{1}\left( T_{Z_{0}}\right)
\end{equation*}
pairs to zero with $H^{0}\left( N_{Y_{0}\backslash Z_{0}}\right) $
. Thus if $Y_{0}$ (infinitesimally) parametrizes some subspace $%
H^{1}\left( \Omega _{Z_{0}}^{2}\right) $, then $Y_{0}$ deforms
(infinitesimally) over $Z'$ in directions ``normal'' to that
subspace.

If our higher-order operators $\left. \left( \left( p_{Y'}\right) _{*}D_{Y/Y^{\prime \prime }}^{Z/Z'}\right)
\right| _{y'}$ pair nicely with relative cohomology of $%
Z/Z'$, we should expect that higher-order obstruction deformation
of $Y/Y'$ in $Z/Z'$ are computed from the higher
derivatives of some generalized Abel-Jacobi mapping or normal function, at
least when there is ``sufficient'' ambient cohomology around. This is indeed
the case as we show in \S\ref{sec.obst}.

\newpage

\section{Review of $\Dd$-modules}\label{sec.dmods}

Let Y and X be complex analytic varieties and let
\begin{equation*}
\begin{diagram}
\node{Y} \arrow{e,t}{g} \node{X}
\end{diagram}
\end{equation*}
be a holomorphic mapping. Let $\Ee$ be a coherent sheaf on $X$ and $%
\Ff$ a coherent sheaf on $Y$.\bigskip

\begin{dfn}\label{dfn.11}
A (left) differential operator on $\Ee$ with
values in $\Ff$ is a $\Bbb{C}$-linear map
\begin{equation*}
\begin{diagram}
\node{g^{-1}\Ee} \arrow{e,t}{T} \node{\Ff}
\end{diagram}
\end{equation*}
such that, for some $n\ge 0$ and for arbitrary sections $f_{0},\ldots ,f_{n}$
of $\Oo_{X}$, the $\left( n+1\right) $-fold commutator
\begin{equation*}
\left[ \ldots \left[ \left[ T,f_{0}\right] ,f_{1}\right] \ldots ,f_{n}\right]
\end{equation*}
is zero, where
\begin{equation*}
\left[ T,f\right](e) = T\left( fe\right) -fT\left( e\right) .
\end{equation*}
The minimum such $n$ is called the \textit{order} of the differential
operator $T$, and we denote the set of all such operators as
\begin{equation*}
\Dd\left( \Ee_{X},\Ff_{Y}\right) .
\end{equation*}
\end{dfn}

Let $\Dd\left( \Ee_{X}\right) $ denote $\Dd\left( \Ee%
_{X},\Ee_{X}\right) $ for the identity map on $X$. Given $%
h:Z\rightarrow Y$ and a coherent sheaf $\Gg$ on $Z$, we have a natural
composition
\begin{equation*}
\Dd\left( \Ff_{Y},\Gg_{Z}\right) \times h^{-1}\Dd\left(
\Ee_{X},\Ff_{Y}\right) \arrow{e} \Dd\left( \Ee_{X},%
\Gg_{Z}\right)
\end{equation*}
so that $\Dd\left( \Ee_{X},\Ff_{Y}\right) $ has the natural
structure of a right $g^{-1}\Dd\left( \Ee_{X}\right) $-module
and a left $\Dd\left( \Ff_{Y}\right) $-module. (We will denote
compositions of operators as $ST$, which means ``do $T$ first, then $S$.'')

For a complex manifold $X$ of dimension $n$, let $\Dd_{X}=\Dd%
\left( \Oo_{X},\Oo_{X}\right) $. The canonical bundle $\omega
_{X}$ of $X$ is a right $\Dd_{X}$-module as follows. Sections of $%
\Ll_{X}$ act on $\omega _{X}$ by Lie differentiation. Multiplying this
action by $\left( -1\right) $, we extend to a right action of the
enveloping algebra $\Uu_{X}$ associated to the Lie algebra $\Ll_{X}$ on $%
\omega _{X}$. But this enveloping algebra, considered as a ring of
(left) differential operators on $\Oo_{X}$, maps to $\Dd_{X}$.
One checks by direct computation in local coordinates that the kernel of
this map acts trivially (on the right) on $\omega _{X}$.

This last construction lets one pass from left $\Dd_{X}$-modules to
right $\Dd_{X}$-modules and back again. One defines this
correspondence on $\Ll_{X}$-modules as follows:

Left to right:
\begin{eqnarray*}
\Ee &\rightarrow& \Ee\tensor _{\Oo_{X}}\omega _{X} \\
\left( f\tensor \omega \right) \cdot \xi &:=& -\xi \left( f\right) \tensor \omega +f\tensor \left( \omega \cdot \xi \right) 
\end{eqnarray*}

Right to left:
\begin{eqnarray*}
\Ff &\rightarrow& \sheafHom_{\Oo_{X}}\left( \omega _{X},\Ff\right) \\
\left( \xi \left( \lambda \right) \right) \left( \omega \right)
	&:=& \lambda \left( \omega \cdot \xi \right) -\lambda \left( \omega \right) \cdot \xi 
\end{eqnarray*}
and extend the action as above to $\Uu_{X}\rightarrow \Dd_{X}$.
(See~\cite[VI,3.3]{bib.B}.)

One relates the structure sheaf $\Oo_{X}$ of holomorphic functions to
the constant sheaf $\Bbb{C}_{X}$ on $X$ by viewing $\Bbb{C}_{X}$
as the ``sheaf of solutions'' associated to the (left) $\Dd_{X}$-module
$\Oo_{X}$, that is, the set of section of $\Oo_{X}$
annihilated by the distinguished ideal of $\Dd_{X}$ generated by the
sheaf $\Ll_{X}$ of vector fields on $X$. The corresponding functor from the
category of bounded $\Dd_{X}$-modules to the bounded derived category
of sheaves is called the ``de Rham functor'' (see~\cite{bib.B}). In the case of $%
\Oo_{X}$ we make the construction explicitly below.\bigskip

\begin{lem}\label{lem.12}
Via the natural map
\begin{equation*}
\Dd\left( \Oo_{X},\Omega _{X}^{*}\right) =\sheafHom\left( \Oo%
_{X},\Omega _{X}^{*}\right) \tensor _{\Oo_{X}}\Dd_{X}=\Omega
_{X}^{*}\tensor _{\Oo_{X}}\Dd_{X}\arrow{e} \omega
_{X}\tensor _{\Oo_{X}}\Dd_{X}\arrow{e} \omega _{X}
\end{equation*}
the sequence
\begin{equation*}
0 \arrow{e} \Dd\left( \Oo_{X},\Oo_{X}\right) \arrow{e,t}{\partial \circ }
\cdots \arrow{e,t}{\partial \circ} \Dd\left( \Oo_{X},\Omega _{X}^{n-1}\right) \arrow{e,t}{\partial\circ} \Dd\left( \Oo_{X},\omega _{X}\right) \arrow{e} \omega _{X} \arrow{e} 0
\end{equation*}
becomes a locally free resolution of $\omega _{X}$ as a right $\Dd_{X}$-module.
In the same way, via
\begin{equation*}
\Dd\left( \omega _{X},\omega _{X}\right) =\Dd%
_{X}\tensor _{\Oo_{X}}\sheafHom\left( \omega
_{X},\omega _{X}\right) =\Dd_{X}\tensor _{\Oo_{X}}%
\Oo_{X}\arrow{e} \Oo_{X}
\end{equation*}
we have
\begin{equation*}
0 \arrow{e} \Dd\left( \omega _{X},\Oo_{X}\right) \arrow{e,t}{\partial \circ }
\cdots \arrow{e,t}{\partial \circ} \Dd\left( \omega _{X},\Omega _{X}^{n-1}\right) \arrow{e,t}{\partial \circ} \Dd\left( \omega _{X},\omega _{X}\right) \arrow{e} \Oo_{X} \arrow{e} 0
\end{equation*}
becomes a locally free resolution of $\Oo_{X}$ as a left $\Dd_{X}$-module.
\end{lem}

\textbf{Proof:} Filter negatively by degree of forms and positively by
order of operator and use the Koszul lemma from elementary homological
algebra, or see~\cite[VI,3.5]{bib.B}. For the second case use that Lie
differentiation commutes with exterior differentiation, and define the left
action of a vector field $\xi $ on $\lambda \in \Dd\left( \omega
_{X},\Omega _{X}^{p}\right) $ by
\begin{equation*}
\left( \lambda \cdot \xi \right) \left( \eta \right) =\lambda \left( \eta
\cdot \xi \right) .
\end{equation*}
\bigskip

One therefore obtains the action of the deRham functor $\omega _{X}%
\overset{L}{\tensor }_{\Dd_{X}}$ on $\Oo_{X}$ by acting on the
above resolution of $\omega _{X}$ with $\tensor _{\Dd_{X}}\Oo
_{X}$. But this computation gives back the deRham complex
\begin{equation*}
0\arrow{e} \Oo_{X}\arrow{e,t}{\partial }\Omega _{X}^{1}%
\arrow{e,t}{\partial } \cdots \arrow{e,t}{\partial} \omega _{X}\arrow{e} 0.
\end{equation*}
So, by the Poincar\'{e} lemma,
\begin{equation*}
\mathcal{H}\left( \omega _{X}\overset{L}{\tensor }_{\Dd_{X}}\Oo_{X}\right) =\Bbb{C}_{X}.
\end{equation*}

If
\begin{equation*}
\begin{diagram}
\node{W} \arrow{e,t}{i} \node{X}
\end{diagram}
\end{equation*}
is the inclusion of $W$ as a closed submanifold of $X$, then $\Dd\left(
\Oo_{X},\Oo_{W}\right) $ is a left $\Dd_{W}$-module
and a right $\Dd_{X}$-module. In the (standard) notation
of~\cite{bib.B}, $\Dd\left( \Oo_{X},\Oo_{W}\right) =\Dd%
_{W\rightarrow X}$. One reverses the sides on which the operators act
by using the left-switching described above, thereby constructing a left $%
\Dd_{X}$-module and right $\Dd_{W}$-module:
\begin{eqnarray*}
\Dd_{X\leftarrow W}
	&=& \sheafHom_{\Oo_{X}}\left( \omega_{X},\Dd\left( \Oo_{X},\Oo_{W}\right) \right)
		\tensor _{\Oo_{W}}\omega _{W} \\
	&=& \Dd\left( \omega _{X},\Oo_{W}\right) \tensor _{\Oo _{W}}\omega _{W} \\
	&=& \Dd\left( \omega _{X},\omega _{W}\right) 
\end{eqnarray*}

Now let $\Ee$ be a left $\Dd_{X}$-module which is a locally
free $\Oo_{X}$-module of finite rank. By Kashiwara's
theorem (\cite[VI, Theorem 7.13]{bib.B})
we can compute the local cohomology sheaf $R\mathcal{H}%
_{W}^{0}\left( \Ee\right) $:
\begin{eqnarray*}
\left( i_{+}\circ i^{!}\right) \Ee
	&=& \Dd_{X\leftarrow W}\overset{L}{\tensor }_{\Dd_{W}}\Dd_{W\rightarrow X}\overset{%
		L}{\tensor }_{\Dd_{X}}\Ee \\
	&=& \left( \Dd\left( \omega _{X},\omega _{W}\right) \overset{L}{%
		\tensor }_{\Dd_{W}}\Dd\left( \Oo_{X},\Oo
		_{W}\right) \right) \overset{L}{\tensor }_{\Dd_{X}}\Ee \\
	&=& \left( \Dd\left( \omega _{X},\omega _{W}\right) \overset{L}{%
		\tensor }_{\Dd_{W}}\left( \Oo_{W}\tensor _{\Oo_{X}}\Dd\left(
		\Oo_{X},\Oo_{X}\right) \right) \right) \overset{%
		L}{\tensor }_{\Dd_{X}}\Ee \\
	&=& \left( \Dd\left( \omega _{X},\omega _{W}\right) \overset{L}{%
		\tensor }_{\Dd_{W}}\Oo_{W}\right) \tensor _{\Oo_{X}}%
		\Dd\left( \Oo_{X},\Oo_{X}\right) \overset{L}{
		\tensor }_{\Dd_{X}}\Ee \\
	&=& \left( \Dd\left( \omega _{X},\omega _{W}\right) \overset{L}{%
		\tensor }_{\Dd_{W}}\Oo_{W}\right) \tensor _{\Oo_{X}}\Ee
		\\
	&\underset{Lemma~\ref{lem.12}}{=}& \left( \Dd\left( \omega _{X},\omega
		_{W}\right) \tensor _{\Dd_{W}}\Dd\left( \omega
		_{W},\Omega _{W}^{*}\right) \right) \tensor _{\Oo_{X}}\Ee \\
	&=& \Dd\left( \omega _{X},\Omega _{W}^{*}\right) \tensor _{\Oo_{X}}\Ee \\
	&=& \Dd\left( \omega _{X}\tensor _{\Oo_{X}}\Ee^{\vee},\Omega _{W}^{*}\right) 
\end{eqnarray*}
(Here `` $^{\vee }$ '' denotes the dual as an $\Oo_{X}$-module. To
check twists, let $W=X$. Via Lemma~\ref{lem.12} we obtain the tautology
\begin{equation*}
\Omega _{X}^{p}=\Dd\left( \Omega _{X}^{n-p},\Omega _{W}^{*}\right) =%
\sheafHom\left( \Omega _{X}^{n-p},\omega _{X}\right)
\end{equation*}
in that case.) Hence in the case of the left $\Dd_{X}$-module $%
\Ee=\Oo_{X}$ we obtain the isomorphism
\begin{equation*}
R\mathcal{H}_{W}^{0}\left( \Oo_{X}\right) =\Dd\left( \omega
_{X},\Omega _{W}^{*}\right)
\end{equation*}
of right $\Dd\left( \omega _{X}\right) $-modules. Thus we have:
\begin{eqnarray*}
R\mathcal{H}_{W}^{0}\left( \Omega _{X}^{p}\right)
	&=& R\mathcal{H}_{W}^{0}\left( \Oo_{X}\right) \tensor _{\Oo_{X}}\Omega_{X}^{p} \\
	&=& \Dd\left( \omega _{X},\Omega _{W}^{*}\right) \tensor _{\Oo_{X}}\Omega _{X}^{p} \\
	&=& \Dd\left( \left( \wedge^{p}T_{X}\right) \tensor _{\Oo_{X}}\omega _{X},\Omega_{W}^{*}\right) \\
	&=& \Dd\left( \Omega _{X}^{n-p},\Omega _{W}^{*}\right) 
\end{eqnarray*}
We will need to check that the isomorphism~(\ref{eq.13}) commutes with the right
action induced by exterior differentiation $\circ \partial _{X}$, that
is, that the diagram
\begin{equation*}
\begin{diagram}
\node{R\mathcal{H}_{W}^{0}\left( \Omega _{X}^{p}\right) } \arrow{s,r}{ \circ \partial } \arrow{e,=} \node{\Dd\left( \Omega _{X}^{n-p},\Omega _{W}^{*}\right) } \arrow{s,r}{ \circ \partial } \\
\node{R\mathcal{H}_{W}^{0}\left( \Omega _{X}^{p+1}\right) } \arrow{e,=} \node{\Dd\left( \Omega _{X}^{n-p-1},\Omega _{W}^{*}\right) } 
\end{diagram}
\end{equation*}
is commutative. Since this is a local computation we may fix a set of local
coordinates $\left( x_{1},\ldots ,x_{n}\right) $. Putting
\begin{equation*}
dx_{1}\wedge \ldots \wedge dx_{n}\leftrightarrow 1
\end{equation*}
equates $\Oo_{X}$ with $\omega _{X}$ as right $\Dd\left(
\omega _{X}\right) = \Dd_{X}$-modules
(with the adjoint action)
and putting
\begin{equation*}
dx_{K}\leftrightarrow \left\{ \left( \sgn\left( K,K^{c}\right)
\right) dx_{K^{c}}\longmapsto 1\right\}
\end{equation*}
equates $\Omega _{X}^{p}$ and $\sheafHom\left( \Omega _{X}^{n-p},\Oo%
_{X}\right) $ as left $\Oo_{X}$-modules. We then compute:
\begin{eqnarray*}
\partial \left( fdx_{K}\right)
	&=& \sum\nolimits_{k\in K^{c}}\frac{\partial f}{\partial x_{k}}dx_{k}\wedge dx_{K} \\
	&\leftrightarrow& \sum\nolimits_{k\in K^{c}}\left( \sgn\left( \left\{ k,K\right\} ,\left\{
		k,K\right\} ^{c}\right) \right) \frac{\partial f}{\partial x_{k}}%
		dx_{\left\{ k,K\right\} ^{c}} \\
	&=& \left\{ \sgn\left( K,K^{c}\right) dx_{K^{c}}\longmapsto f\right\} \circ \partial \\
	&=& \sgn\left( K,K^{c}\right) \sum\nolimits_{k\in K^{c}}\left\{
		dx_{K^{c}-\left\{ k\right\} }\longmapsto \frac{\partial f}{\partial x_{k}}\right\} 
\end{eqnarray*}
Thus
\begin{equation}
\left( R\mathcal{H}_{W}^{0}\left( \Omega _{X}^{*}\right) ,\partial \right)
=\left( \Dd\left( \Omega _{X}^{n-*},\Omega _{W}^{*}\right) ,\circ
\partial \right) .
\label{eq.13}
\end{equation}

Formula~(\ref{eq.13}) is a bit removed from the classical treatment of local
cohomology in algebraic geometry, but a local, non-canonical correspondence
is not hard to see. Namely, if we have $W$ realized as a complete
intersection of hypersurfaces
\begin{equation*}
W=Y_{1}\cap \;\ldots \;\cap Y_{c}
\end{equation*}
then for the locally free $\Oo_{X}$-module $\Omega _{X}^{p}$,
local cohomology is computed from the sheaf
\begin{equation*}
\frac{\Omega _{X}^{p}\left( \infty Y_{1}+\ldots +\infty
Y_{c}\right) }{\sum\nolimits_{j}\Omega _{X}^{p}\left( \infty
Y_{1}+\ldots +\infty Y_{j-1}+\infty Y_{j+1}+\ldots +\infty
Y_{c}\right) }.
\end{equation*}
On the other hand, if we fix a holomorphic projection
\begin{equation*}
\begin{diagram}
\node{X} \arrow{e,t}{p} \node{W}
\end{diagram}
\end{equation*}
as left inverse of the inclusion $i$, there is a well-defined residue
operator
\begin{equation*}
\begin{diagram}
\node{\frac{\omega _{X}\left( \infty Y_{1}+\ldots +\infty
Y_{c}\right) }{\sum\nolimits_{j}\omega _{X}\left( \infty
Y_{1}+\ldots +\infty Y_{j-1}+\infty Y_{j+1}+\ldots +\infty
Y_{c}\right) }}
\arrow{e,t}{\text{res}}
\node{\omega _{W}}
\end{diagram}
\end{equation*}
and so a well-defined map
\begin{eqnarray*}
\frac{\Omega _{X}^{p}\left( \infty Y_{1}+\ldots +\infty Y_{c}\right)
		}{\sum\nolimits_{j}\Omega _{X}^{p}\left( \infty Y_{1}+\ldots
		+\infty Y_{j-1}+\infty Y_{j+1}+\ldots +\infty Y_{c}\right) }
	&\longrightarrow& \Dd\left( \Omega _{X}^{n-p},\omega _{W}\right) \\
\eta
	&\longmapsto& \left( \psi \mapsto res\left( \psi \wedge \eta \right) \right) 
\end{eqnarray*}
which induces~(\ref{eq.13}).

\newpage

\section{The sheaf of deRham normal differential operators}\label{sec.ndo}

\subsection{Definition of $\Dd\left( \Omega _{Z/Z'},\Omega_{Y/Y'}\right) $}\label{sec.2.1}
Throughout this paper, we assume that we have a commutative diagram
\begin{equation}
\begin{diagram}
\node{Y} \arrow{s,r}{ p_{Y'}} \arrow{e,t}{h} \node{Z} \arrow{s,r}{ p_{Z'}} \\
\node{Y'} \arrow{e,t}{h'} \node{Z'} 
\end{diagram}
\label{eq.211}
\end{equation}
of holomorphic mappings of complex analytic manifolds. We assume that
vertical maps are (surjective) smooth projections with connected
fibers.\bigskip

\begin{dfn}\label{dfn.212}
We call a diagram~(\ref{eq.211}) relatively immersive
if
\begin{equation*}
Y\arrow{e} Y'\times _{Z'}Z
\end{equation*}
is a closed immersion.\bigskip
\end{dfn}

So we have a fibered product:
\begin{equation*}
\begin{diagram}
\node{X:=X'\times _{Z'}Z} \arrow{s,r}{ p_{X'}} \arrow{e,t}{j} \node{Z} \arrow{s,r}{ p_{Z'}} \\
\node{X'=Y'\times Z'} \arrow{e,t}{j'} \node{Z'} 
\end{diagram}
\end{equation*}
We then have:
\begin{equation*}
\begin{diagram}
\node{Y} \arrow{s,r}{ p_{Y'}} \arrow{e,t}{i} \node{X} \arrow{s,r}{ p_{X'}} \arrow{e,t}{j} \node{Z} \arrow{s,r}{ p_{Z'}} \\
\node{Y'} \arrow{e,t}{i'} \node{X'} \arrow{e,t}{j'} \node{Z'} 
\end{diagram}
\end{equation*}
If~(\ref{eq.211}) is relatively immersive, then $i$ is a closed immersion and $j$
is a smooth projection.

Let:
\begin{eqnarray*}
n &=& \mathop{\mathrm{rel}}\dim Z/Z' \\
d &=& \mathop{\mathrm{rel}}\dim Y/Y' \\
c &=& n-d 
\end{eqnarray*}
We will work throughout with sheaves of relative K\"{a}hler differentials.
For example, on $Z$ these are given by the formulas:
\begin{eqnarray*}
\Omega _{Z/Z'}^{0} &=& \Oo_{Z} \\
\Omega _{Z/Z'}^{q} &=& \frac{\Omega _{Z}^{q}}{p_{Z^{\prime }}^{*}\left(
	\Omega _{Z'}^{1}\right) \wedge \Omega _{Z}^{q-1}}\quad q>0 \\
\end{eqnarray*}
We want to apply our discussion of sheaves of differential operators in \S\ref{sec.dmods}
to the case
\begin{equation*}
\Ee_{X}=\Omega _{X/X'}^{*}=j^{*}\Omega _{Z/Z'}^{*}\hspace{0.5in}\Ff_{Y}=\Omega _{Y/Y'}^{*}.
\end{equation*}

We define the bicomplex
\begin{equation*}
\Dd\left( \Omega _{Z/Z'},\Omega _{Y/Y'}\right) ^{p,q}=\left\{ \left( \lambda ^{p,q}\in
\Dd\left( \Omega_{X/X'}^{n-p},\Omega _{Y/Y'}^{q}\right) \right)
: \text{$\lambda^{p,q}$ is $\Oo_{Y}$-linear}\right\}
\end{equation*}
with total differential
\begin{equation*}
\left( \circ \partial _{Z/Z'}\right) +\left( -1\right) ^{p}\left(
\partial _{Y/Y'}\circ \right)
\end{equation*}
made up of two summands of degree $\left( 1,0\right) $ and $\left(
0,1\right) $ respectively.\smallskip

\subsection{Structure of $\Dd\left( \Omega _{Z/Z'},\Omega_{Y/Y'}\right) $}\label{sec.2.2}

As in \S\ref{sec.dmods}, for any complex manifold $W$, let
\begin{equation*}
\Ll_{W}
\end{equation*}
denote the sheaf of Lie algebras of holomorphic vector fields on $W$. The
product structure on $X'=Y'\times Z'$ above
distinguishes a sub-Lie-algebra of $\Ll_{X'}$, namely the subsheaf of
the vector fields tangent to the fibers of the projection to $Y'$.
We denote this subsheaf of Lie algebras as
\begin{equation*}
\Ll_{Y',Z'}.
\end{equation*}
We define
\begin{equation}
\Ll_{Y',Z}
\label{eq.221}
\end{equation}
as the sub-Lie-algebra of $\Ll_{X}$ made up of elements which are liftings of
elements of $\Ll_{Y',Z'}$, that is, vector fields $\xi $
for which there exists a vector field $\xi '$ in $\Ll_{Y',Z'}$ such that
\begin{equation*}
\xi \left( g\circ p_{X}\right) =\xi \left( g\right) \hspace{0.5in}g\in
\Oo_{X}.
\end{equation*}
Then $\Ll_{Y',Z}$ acts on the right on the complex $\Dd\left( \Omega
_{Z/Z'},\Omega _{Y/Y'}\right) $. $\Ll_{Y',Z}$ has a subsheaf
\begin{equation*}
\Ll_{Y',Z}^{0}
\end{equation*}
of vector fields which are lifts of the zero vector field, that is, vector
fields pointing along the fibers of $p_{X'}$. Notice that $%
\Ll_{Y',Z}$ is an $\Oo_{X'}$-module but that the
submodule $\Ll_{Y',Z}^{0}$ is (compatibly) an $\Oo_{X}$-module.
We will abbreviate these two facts by referring to the `` $\Oo%
_{X/X'}$-module structure'' of $\Ll_{Y',Z}$. Notice that
\begin{equation}
\left[ \Ll_{Y',Z}^{0},\Ll_{Y',Z}\right] \subseteq
\Ll_{Y',Z}^{0}
\label{eq.222}
\end{equation}
and, near $x\in X$,
\begin{equation*}
\left( \Ll_{Y',Z}\right) _{x}\cong \left( \Ll_{Y',Z}^{0}\right) _{x}\oplus p_{X'}^{-1}\left( j'\left( \Ll_{Z'}\right) \right) _{x}.
\end{equation*}

Lie and exterior differentiation commute. Thus the complex $\Dd\left(
\Omega _{Z/Z'},\Omega _{Y/Y'}\right) ^{*,*}$ is
a right $\Ll_{Y',Z}$-module. In fact, let
\begin{equation*}
\text{$\Uu$, resp.~$\Uu^{0}$,}
\end{equation*}
denote the $p_{X'}^{-1}\Oo_{X'}$-enveloping
algebras associated to the Lie algebra $\Ll_{Y',Z}$, respectively $%
\Ll_{Y',Z}^{0}$ . Similary let
\begin{equation*}
\Uu'
\end{equation*}
denote the $\Oo_{X'}$-enveloping algebra associated to $%
\Ll_{Y',Z'}$. Then the complex $\Dd\left( \Omega
_{Z/Z'},\Omega _{Y/Y'}\right) ^{*,*}$ becomes a right $\Uu$-module.\smallskip

\subsection{Structure of $\Dd\left( \Omega _{Z/Z'},\Omega_{Y/Y'}\right) _{\Jj}$}\label{sec.2.3}

If $\Jj'$ is a sheaf of ideals on $X'$ which
annihilates $i'\left( Y'\right) $ and we define
\begin{equation*}
\Jj=p_{X'}^{*}\Jj',
\end{equation*}
then we can define the complex:
\begin{equation*}
\Dd\left( \Omega _{Z/Z'},\Omega _{Y/Y'}\right) _{\Jj}^{p,q}=\left\{ \lambda \in \Dd\left( \Omega
_{Z/Z'},\Omega _{Y/Y'}\right) ^{p,q}:\lambda
\left( g\eta \right) =0,\forall g\in \Jj,\eta \in \Omega _{X/X'}^{n-p}\right\}
\end{equation*}
Also let
\begin{equation}
\begin{array}{rcl}
\Omega _{X_{\Jj}/X_{\Jj'}'}^{*} &=& \frac{\Omega _{X/X'}^{*}}{p_{X}^{*}\Jj\cdot \Omega _{X/X'}^{*}} \\
\Uu_{\Jj'}' &=& \left\{ u'\in \Uu^{\prime }:u'\left( g\right) |_{Y'}=0,\forall g\in \Jj'\right\} \\
\Uu_{\Jj} &=& \left\{ u\in \Uu:u\left( \eta \right) |_{Y}=0,\eta \in \Jj\cdot \Omega _{X/X'}^{*}\right\} 
\end{array}
\label{eq.233}
\end{equation}
where ``$|_{Y}$'' denotes the pullback $i^{*}$. Let $\rho _{Y}$
denote the pullback via $h$ composed with restriction of forms to $Y$. Then
\begin{equation*}
\rho _{Y}\cdot \Uu_{\Jj}\subseteq \Dd\left( \Omega _{X_{\Jj%
}/X_{\Jj'}},\Omega _{Y/Y'}\right) ^{p,n-p}=%
\Dd\left( \Omega _{X/X'},\Omega _{Y/Y'}\right) _{\Jj}^{p,n-p},\forall p.
\end{equation*}
\smallskip

\subsection{Structure of $\Dd^{\bot }\left( \Omega _{Z/Z'}^{p},\omega _{Y/Y'}\right) $}\label{sec.2.4}

Next we let
\begin{equation*}
\Dd_{Y/Y'}
\end{equation*}
be the algebra of $\Oo_{Y'}$-linear differential operators on $%
Y$. As we saw in \S\ref{sec.dmods}, $\Dd_{Y/Y'}$ acts on the right on the
relative dualizing sheaf
$\omega _{Y/Y'}^{} $ and we can form
\begin{equation}
\Dd^{\bot }\left( \Omega _{Z/Z'}^{p},\omega _{Y/Y'}\right) =\omega _{Y/Y'}\tensor _{\Dd_{Y/Y'}}%
\Dd\left( \Omega _{Z/Z'}^{p},\Oo_{Y}\right) .
\label{eq.241}
\end{equation}

To compute the derived functors of~(\ref{eq.241}), let $\Dd_{Y/Y'}\left( \Oo_{Y},\Omega _{Y/Y'}^{*}\right) $ denote the sheaf
of $\Oo_{Y'}$-linear elements of $\Dd\left( \Oo%
_{Y},\Omega _{Y/Y'}^{*}\right) $ and form:
\begin{equation*}
0\arrow{e} \Dd_{Y/Y'}\left( \Oo_{Y}\right) \arrow{e,t}{\partial\circ} \ldots \arrow{e,t}{\partial\circ}\Dd_{Y/Y'}\left( \Oo_{Y},\Omega
_{Y/Y'}^{d-1}\right) \arrow{e,t}{\partial \circ }\Dd%
_{Y/Y'}\left( \Oo_{Y},\omega _{Y/Y'}\right) \arrow{e}
\omega _{Y/Y'}\arrow{e} 0
\end{equation*}
By Lemma~\ref{lem.12}, this is a (locally) free resolution of $\omega _{Y/Y'}$ as a right $\Dd_{Y/Y'}$-module, so we have the
equality
\begin{equation}
\begin{array}{rcl}
\omega _{Y/Y'}\overset{L}{\tensor }_{\Dd_{Y/Y'}}%
\Dd\left( \Omega _{Z/Z'}^{n-p},\Oo_{Y/Y'}\right)
&=& \mathcal{H}\left( \Dd\left( \Oo_{Y/Y'},\Omega
_{Y/Y'}^{*}\right) \tensor _{\Dd_{Y/Y'}}\Dd%
\left( \Omega _{Z/Z'}^{n-p},\Oo_{Y/Y'}\right) \right) \\
&=& \mathcal{H}\left( \Dd\left( \Omega _{Z/Z'},\Omega
_{Y/Y'}\right) ^{p,*}\right) 
\end{array}
\label{eq.242}
\end{equation}
of right $\Uu$-modules. Similarly we let
\begin{equation*}
\Dd_{Z/Z'}\left( \Omega _{Z/Z'}^{n-*},\Omega
_{Y/Y'}^{*}\right)
\end{equation*}
denote the $\Oo_{Z'}$-linear elements of $\Dd\left(
\Omega _{Z/Z'}^{n-*},\Omega _{Y/Y'}^{*}\right) $ and form
the right $\Uu^{0}$-module
\begin{equation*}
\Dd_{Z/Z'}^{\bot }\left( \Omega _{Z/Z'}^{n-*},\omega _{Y/Y'}\right)
\end{equation*}
as above.\smallskip

\subsection{Cohomology of $\Dd\left( \Omega _{Z/Z'},\Omega _{Y/Y'}\right) ^{*,*}$}\label{sec.2.5}

In this section we assume that
\begin{equation*}
\begin{diagram}
\node{Y} \arrow{s,r}{ p_{Y'}} \arrow{e,t}{h} \node{Z} \arrow{s,r}{ p_{Z'}} \\
\node{Y'} \arrow{e,t}{h'} \node{Z'} 
\end{diagram}
\end{equation*}
in~(\ref{eq.211}) is relatively immersive. We can then explicitly compute the
cohomology of the complex $\Dd\left( \Omega _{Z/Z'},\Omega
_{Y/Y'}\right) ^{*,*}$. Referring to~(\ref{eq.242}), the $E_{1}$-term in
the ``$p$-filtration'' spectral sequence associated to the bicomplex
\begin{equation*}
\sum\nolimits_{p,q}\Dd\left( \Omega _{Z/Z'},\Omega
_{Y/Y'}\right) ^{p,q}
\end{equation*}
is exactly
\begin{equation*}
\omega _{Y/Y'}\overset{L}{\tensor }_{\Dd_{Y/Y'}}%
\Dd\left( \Omega _{Z/Z'}^{n-*},\Oo_{Y/Y'}\right)
.
\end{equation*}
Since~(\ref{eq.211}) is relatively immersive, $\Dd\left( \Omega _{Z/Z'}^{n-*},\Oo_{Y/Y'}\right) $ is a locally free $\Dd_{Y/Y'}$-module:
\begin{equation*}
\begin{tabular}{l}
$E_{1}^{p,d}=\Dd^{\bot }\left( \Omega _{Z/Z'}^{n-p},\omega _{Y/Y'}\right) $ \\
$E_{1}^{p,q}=0,\quad q\neq d$
\end{tabular}
\end{equation*}
For the immersion
\begin{equation*}
\begin{diagram}
\node{Y} \arrow{e,t}{i} \node{Y'\times Z}
\end{diagram}
\end{equation*}
we let $\mathcal{H}_{Y\backslash Y'\times Z}^{*}$ denote the
derived functors of the functor which assigns to a sheaf on $X=Y'\times Z$ the subsheaf of the pullback consisting of sections supported on
the submanifold $Y$. Now, by the assumption that $i$ is locally a closed
immersion,
\begin{equation*}
\mathcal{H}_{Y\backslash Y'\times Z}^{k}=0
\end{equation*}
for $k\neq n+c$ . Also, by~(\ref{eq.13}), we have
\begin{equation}
\Dd^{\bot }\left( \Omega _{Z/Z'}^{n-p},\omega
_{Y/Y'}\right) \cong \mathcal{H}_{Y\backslash Y'\times
Z}^{n+c}\left( \Omega _{Z/Z'}^{p}\right)
\label{eq.252}
\end{equation}
We next filter $\Uu$ and $\Dd\left( \Omega _{Z/Z'}^{n-*},\omega
_{Y/Y'}\right) $ respectively by assigning to
\begin{equation*}
\begin{tabular}{l}
$P\in \Uu$ \\
$Q\in \Dd\left( \Omega _{Z/Z'}^{n-*},\omega _{Y/Y^{\prime }}\right) $
\end{tabular}
\end{equation*}
an order $r$ equal to the order of the operator
\begin{equation*}
\left( f\longmapsto P\left( f\right) \right) ,\;f\in \Oo_{Z'},
\end{equation*}
respectively, the maximal order
over all $\omega \in \Omega _{Z/Z'}^{n-*}$ of the operator
\begin{equation*}
\left( f\longmapsto Q\left( f\omega \right) \right) ,\;f\in \Oo_{Z'}.
\end{equation*}
We let
\begin{equation}
\begin{array}{rcl}
_{\left[ r\right] }\Uu
	&=& \mathop{\mathrm{image}}\,\left\{ P\in \Uu:\deg P\le r\right\} \\
_{\left[ r\right] }\Dd^{\bot }\left( \Omega _{Z/Z'}^{n-*},\omega _{Y/Y'}\right)
	&=&\mathop{\mathrm{image}}\,\left\{ P\in \Dd\left(
		\Omega _{Z/Z'}^{n-*},\omega _{Y/Y'}\right) :\deg P\le r\right\}. 
\end{array}
\label{eq.245}
\end{equation}
The $r$-th graded quotient of the respective filtrations are
\begin{equation}
\begin{array}{c}
\Uu^{0}\tensor S^{r}T_{Z'} \\
\Dd_{Z/Z'}^{\bot }\left( \Omega _{Z/Z'}^{n-*},\omega
_{Y/Y'}\right) \tensor S^{r}T_{Z'}
\end{array}
\label{eq.253}
\end{equation}
where $S^{r}$ denotes the $r$-th symmetric product. Thus we have a
filtration:
\begin{equation}
E_{2}^{pq}=\mathcal{H}^{p}\left( \Dd^{\bot }\left( \Omega _{Z/Z'}^{n-*},\omega _{Y/Y'}\right) ,\circ \partial \right) \supseteq \ldots \supseteq \mathcal{H}_{r}^{p}\supseteq \ldots \supseteq \mathcal{H}_{0}^{p}
\label{eq.254}
\end{equation}
where
\begin{eqnarray*}
\frac{\mathcal{H}_{r}^{p}}{\mathcal{H}_{r-1}^{p}}
	&=& \mathcal{H}^{p}\left(\Dd_{Z/Z'}^{\bot }\left( \Omega _{Z/Z'}^{n-*},\omega
		_{Y/Y'}\right) \tensor S^{r}T_{Z'},\circ \partial \right) \\
	&\cong& \mathcal{H}_{Y\backslash Y'\times _{Z'}Z}^{c+p}\left( p_{Z'}^{-1}S^{r}T_{Z'}\right) \\
	&=& \left\{ _{\text{$i^{*}S^{r}T_{Z'}$ if $p=c$}}^{\text{0 if $p\neq c$}}\right. 
\end{eqnarray*}
\smallskip

\subsection{Definition of $\Dd_{Y/Y'}^{Z/Z'}$}\label{sec.2.6}

Returning now to the general case~(\ref{eq.211}), we truncate our original complex
\begin{eqnarray*}
\Dd\left( \Omega _{Z/Z'}^{\le d},\Omega _{Y/Y^{\prime }}\right) ^{p,q}
	&=& \Dd\left( \Omega _{Z/Z'}^{n-p},\Omega _{Y/Y'}^{q}\right) \text{{} for $p\ge c$} \\
\Dd\left( \Omega _{Z/Z'}^{\le d},\Omega _{Y/Y^{\prime }}\right) ^{p,q}
	&=& 0 \text{{} for $p<c$} 
\end{eqnarray*}
and we define our main object of study in this paper
\begin{equation}
\Dd_{Y/Y'}^{Z/Z'}=\mathcal{H}^{n}\left( \Dd%
\left( \Omega _{Z/Z'}^{\le d},\Omega _{Y/Y'}\right)
^{*,*},\left( \circ \partial \right) +\left( -1\right) ^{p}\left( \partial
\circ \right) \right) .
\label{eq.261}
\end{equation}
$\Dd_{Y/Y'}^{Z/Z'}$ is clearly a right $\Uu$-module.
Notice that the map
\begin{equation*}
\begin{diagram}
\node{\Omega _{Z/Z'}^{\le d}} \arrow{e,t}{\rho _{Y}} \node{\Omega_{Y/Y'}}
\end{diagram}
\end{equation*}
given by pulling back differential forms is closed but not exact, and that
we have a natural map
\begin{equation*}
\Dd^{\bot }\left( \Omega _{Z/Z'}^{d+1},\omega _{Y/Y'}\right) =\omega _{Y/Y'}\overset{L_{0}}{\tensor }_{\Dd%
_{Y/Y'}}\Dd\left( \Omega _{Z/Z'}^{d+1},\Oo%
_{Y}\right) \arrow{e,t}{\circ \partial }\Dd%
_{Y/Y'}^{Z/Z'}
\end{equation*}
\smallskip

\subsection{$\Dd_{Y/Y'}^{Z/Z'}$ as local cohomology}\label{sec.2.7}

Notice also that the relative deRham complex
\begin{equation*}
0\arrow{e} \Oo_{Z}\arrow{e,t}{\partial }\Omega
_{Z/Z'}^{1}\arrow{e,t}{\partial }\ldots \arrow{e,t}{\partial} \omega _{Z/Z'}\arrow{e} 0
\end{equation*}
is a resolution of $p_{Z'}^{-1}\Oo_{Z'}$. So if~(\ref{eq.211})
is relatively immersive, we have isomorphisms
\begin{equation*}
\mathcal{H}\left( \Dd\left( \Omega _{Z/Z'},\Omega
_{Y/Y'}\right) ^{*,*},\left( \circ \partial \right) +\left(
-1\right) ^{p}\left( \partial \circ \right) \right) =\mathcal{H}%
_{Y/Y'}^{*}\left( \Bbb{C}_{Y'\times Z'}\right) .
\end{equation*}
Thus, by \S\ref{sec.2.5}, we have
\begin{equation*}
\Dd_{Y/Y'}^{Z/Z'}=\ker \left( \Dd^{\bot }\left(
\Omega _{Z/Z'}^{d},\omega _{Y/Y'}\right) \arrow{e,t}{\partial_{d}} \Dd^{\bot }\left( \Omega _{Z/Z'}^{d-1},\omega _{Y/Y'}\right) \right) .
\end{equation*}
And using the increasing filtration $_{\left[ r\right] }\Dd%
_{Y/Y'}^{Z/Z'}$ of $\Dd_{Y/Y'}^{Z/Z'}$ given by~(\ref{eq.245}) and using~(\ref{eq.253}) we have
\begin{equation}
\frac{_{\left[ r\right] }\Dd_{Y/Y'}^{Z/Z'}}{_{\left[
r-1\right] }\Dd_{Y/Y'}^{Z/Z'}}=\mathcal{H}\left(
\Dd_{Z/Z'}\left( \Omega _{Z/Z'}^{\le d},\Omega
_{Y/Y'}\right) ^{*,*}\right) \tensor S^{r}T_{Z'}.
\label{eq.273}
\end{equation}
Also $\Dd_{Y/Y'}^{Z/Z'}$ has a distinguished $\Uu$-submodule given as the image of the mapping
\begin{equation*}
\mathcal{H}\left( \Dd\left( \Omega _{Z/Z'}^{>d}\left[ 1\right]
,\Omega _{Y/Y'}\right) ^{*,*}\right) \arrow{e} \mathcal{H}%
\left( \Dd\left( \Omega _{Z/Z'}^{\le d},\Omega _{Y/Y'}\right) ^{*,*}\right)
\end{equation*}
induced by the morphism of complexes
\begin{equation*}
\Omega _{Z/Z'}^{\le d}\arrow{e,t}{\partial _{d}}%
\Omega _{Z/Z'}^{>d}\left[ 1\right] .
\end{equation*}

In fact~(\ref{eq.254}) tells us that
\begin{equation}
\mathcal{H}\left( \Dd\left( \Omega _{Z/Z'}^{>d}\left[ 1\right]
,\Omega _{Y/Y'}\right) ^{*,*}\right) =\;\Dd^{\bot }\left(
\Omega _{Z/Z'}^{d+1},\omega _{Y/Y'}\right) \circ
\partial _{d}
\label{eq.274}
\end{equation}
and that the sequence
\begin{equation}
0\arrow{e} \Dd^{\bot }\left( \Omega _{Z/Z'}^{d+1},\omega_{Y/Y'}\right) \circ \partial _{d}
\arrow{e} \Dd_{Y/Y'}^{Z/Z'}
\arrow{e} p_{Y'}^{-1} \Dd^{\bot }\left(\Omega _{Z'}^{d'},\omega _{Y'}\right)
\arrow{e} 0
\label{eq.275}
\end{equation}
is exact.
Also
\begin{equation}
_{\left[ 0\right] }\Dd_{Y/Y'}^{Z/Z'}\;=\Oo%
_{Y'}\cdot \rho \oplus \Dd_{Z/Z'}^{\bot }\left(
\Omega _{Z/Z'}^{d+1},\omega _{Y/Y'}\right) \circ
\partial _{d}
\label{eq.276}
\end{equation}
where $\rho $ is the class given by the restriction element $\rho _{Y}$
above.\smallskip

\subsection{Functoriality of $\Dd_{Y/Y'}^{Z/Z'}$}\label{sec.2.8}

Returning to the general case, one easily checks that $\Dd_{Y/Y'}^{Z/Z'}$ satisfies the basic functorial properties:

Any commutative diagram of morphisms
\begin{equation*}
\begin{diagram}
\node{Z_{1}} \arrow{s,r}{ } \arrow{e}  \node{Z_{2}} \arrow{s,r}{ } \\
\node{Z_{1}'} \arrow{e}  \node{Z_{2}'} 
\end{diagram}
\end{equation*}
with the vertical maps smooth and surjective induces a natural map
\begin{equation*}
\Dd_{Y/Y'}^{Z_{1}/Z_{1}'}\arrow{e} \Dd%
_{Y/Y'}^{Z_{2}/Z_{2}'}.
\end{equation*}

Any fibered product
\begin{equation*}
\begin{diagram}
\node{Y_{1}} \arrow{s,r}{ } \arrow{e}  \node{Y_{2}} \arrow{s,r}{ } \\
\node{Y_{1}'} \arrow{e}  \node{Y_{2}'} 
\end{diagram}
\end{equation*}
(with the vertical maps smooth and surjective) induces a natural map
\begin{equation*}
\Dd_{Y_{2}/Y_{2}'}^{Z/Z'}\arrow{e} \Dd%
_{Y_{1}/Y_{1}'}^{Z/Z'}.
\end{equation*}

If $Z=Z'$ and so necessarily $Y=Y'$,
we shall occasionally make the abbreviation
\begin{equation*}
\Dd_{Y}^{Z} := \Dd_{Y/Y}^{Z/Z}.
\end{equation*}
as in equation~(\ref{eq.275}).

Finally, referring to the decomposition~(\ref{eq.276}), the basic formula
\begin{equation*}
L_{\varsigma }=\left\< \zeta \right. \left| \partial \left( {}\right)
\right\> +\partial \left\< \zeta \right. \left| {\ }\right\>
\end{equation*}
for Lie differentiation shows that any element of $\Ll_{Z/Z'}^{0}$
(i.e., the sheaf of vector fields on $Z$ which annihilate functions from $Z'$)
multiplies $\Oo_{Y}\cdot \rho $ into $\Dd_{Z/Z'}^{\bot }\left( \Omega _{Z/Z'}^{d+1},\omega _{Y/Y'}\right) \circ \partial _{d}$.
A more subtle point is that, for a
general element $\xi \in \Ll_{Z/Z'}$,
\begin{equation*}
\rho \cdot \xi =\rho \cdot \left( \partial \left\< \xi \right. \left| {\ }\right\> +\left\< \xi \right. \left| \partial\ \right\>
\right) \notin \Dd^{\bot }\left( \Omega _{Z/Z'}^{d+1},\omega
_{Y/Y'}\right) \circ \partial _{d}.
\end{equation*}
This is because, unless $\xi \in \Ll_{Z/Z'}^{0}$, the component
operations in in the expression
\begin{equation*}
\partial \left\< \xi \right. \left| {\ }\right\> +\left\< \xi \right. \left| \partial \;\right\>
\end{equation*}
cannot be split apart as a sum of well defined operations on relative
differentials. For example $\left\< \xi \right. \left| {\ }\right\> $
is not well-defined on relative differentials unless $\xi \in \Ll_{Z/Z'}^{0}$.

\newpage

\section{Distinguished filtrations, cyclicity}\label{sec.filt}

Throughout this section we continue with the situation of~\S\ref{sec.2.1} and we assume
that $h$ in~(\ref{eq.211}) is a relative immersion. We now consider the standard
increasing filtration
\begin{equation*}
\Uu_{\left[ m\right] }
\end{equation*}
of $\Uu$ given by assigning degree one to elements of $\Ll_{Y',Z}$
(and degree zero to elements of $p_{X'}^{-1}\Oo_{X'}$).
This filtration induces increasing filtration
\begin{equation*}
\left( \Dd_{Y/Y'}^{Z/Z'}\right) _{\left[ m\right] }
\end{equation*}
on
\begin{equation*}
\rho \cdot \Uu\subseteq \Dd_{Y/Y'}^{Z/Z'}.
\end{equation*}
Referring to~(\ref{eq.221}), we can define a ``relative normal bundle'' as the
image of
\begin{equation*}
\frac{i^{*}\Ll_{Y',Z}}{\Ll_{Y',Y}}
\end{equation*}
in $\Dd_{Y/Y'}^{Z/Z'}$ via the map
\begin{equation*}
\zeta \mapsto \rho _{Y}\cdot L_{\zeta }.
\end{equation*}
This image $N_{Y\backslash Z}$ is the second summand of a splitting
\begin{equation*}
\left( \Dd_{Y/Y'}^{Z/Z'}\right) _{\left[ 1\right] }=%
\Oo_{Y}\cdot \rho \cdot\ \oplus \ N_{Y\backslash Z}
\end{equation*}
Notice that $N_{Y\backslash Z}$ is not the usual normal sheaf $%
N_{Y\backslash Y'\times Z}$, but is rather the $\Oo_{Y/Y'}$-module
given by normal vector fields in $X$ which are
liftings of tangent vector fields on $Z$, that is, we have an exact
sequence
\begin{equation*}
0\arrow{e} N_{Y\backslash Y'\times _{Z'}Z}\arrow{e}
N_{Y\backslash Z}\arrow{e} p_{Y'}^{-1}h^{\prime *}T_{Z'}\arrow{e} 0
\end{equation*}
with the sub-object being an $\Oo_{Y}$-module and the quotient an $%
\Oo_{Y'}$-module.\bigskip

\begin{lem}\label{lem.31}
With reference to the filtration $_{\left[ r\right] }%
\Dd_{Y/Y'}^{Z/Z'}$ of $\Dd_{Y/Y'}^{Z/Z'}$ defined in \S\ref{sec.2.7}, $_{\left[ 0\right] }\Dd%
_{Y/Y'}^{Z/Z'}$ is the $\Uu^{0}$-submodule of $\Dd%
_{Z/Z'}^{\bot }\left( \Omega _{Z/Z'}^{d},\omega _{Y/Y'}\right) $ generated by $\rho \cdot\ _{\left[ 0\right] }\Dd%
_{Y/Y'}^{Z/Z'}$ has an increasing filtration induced by
the natural filtration on $\Uu^{0}$. For m $>$ 0, the $m$-th graded quotient
of this filtration is
\begin{equation*}
S_{\Oo_{Y}}^{m}N_{Y\backslash Y'\times _{Z'}Z}.
\end{equation*}
\end{lem}

\textbf{Proof:} Since all computations will be $\Oo_{Y'}$-linear,
we can and will assume throughout this proof that $Y'$ is
a point. We let $\Ll_{Z/Z'}^{0}$ denote the vector fields on $Z$
which annihilate functions from $Z'$ and let $\Uu_{Z/Z'}^{0}$
denote the $\Bbb{C}$-enveloping algebra associated to $\Ll_{Z/Z'}^{0}$%
. We now consider the $\left( \Oo_{Z'},\Ll_{Z/Z'}^{0}\right) $-enveloping algebra
\begin{equation*}
\Uu_{\Oo_{Z}}\left( \Ll_{Z/Z'}^{0}\right)
\end{equation*}
obtained by defining a Lie algebra structure on
\begin{equation*}
\Oo_{Z}+\Ll_{Z/Z'}^{0}
\end{equation*}
by the rule
\begin{equation*}
\left[ f+\lambda ,g+\mu \right] =\left( \lambda \left( g\right) -\mu \left(
f\right) \right) +\left[ \lambda ,\mu \right] ,
\end{equation*}
letting $\Uu\left( \Oo_{Z}+\Ll_{Z/Z'}^{0}\right) $ denote
its universal enveloping algebra over $\Bbb{C}$ and defining
\begin{equation*}
\Uu_{\Oo_{Z}}\left( \Ll_{Z/Z'}^{0}\right) =\frac{\left(
\Oo_{Z}+\Ll_{Z/Z'}^{0}\right) \cdot \Uu\left( \Oo_{Z}+\Ll_{Z/Z'}^{0}\right) }{%
\left\< f\tensor \left( g+\lambda \right) =\left( fg+f\lambda \right)
\right\> }.
\end{equation*}
We consider the standard resolution of $\Oo_{Z}$ as a $\Uu_{\Oo%
_{Z}}\left( \Ll_{Z/Z'}^{0}\right) $-module and the standard
identification of complexes
\begin{equation*}
\sheafHom_{\Uu_{\Oo_{Z}}\left( \Ll_{Z/Z'}^{0}\right)
}\left( \Uu_{\Oo_{Z}}\left( \Ll_{Z/Z'}^{0}\right) \tensor _{%
\Oo_{Z}}\wedge ^{p}\Ll_{Z/Z'}^{0},\ \Oo_{Z}\right) ,D
\end{equation*}
where $D$ is the coboundary induced from the bar resolution. Since
\begin{equation*}
\rho \cdot \partial =0,
\end{equation*}
the standard formulas for Lie differentiation show that we have a morphism
of complexes:
\begin{equation}
\thinD
\begin{diagram}
\node{} \arrow{s} \node{} \arrow{s} \\
\node{\wedge ^{2}\Ll_{Z/Z'}^{0}\tensor _{\Oo_{Z}}\Uu_{\Oo_{Z}}\left( \Ll_{Z/Z'}^{0}\right)} \arrow{s} \arrow{e} \node{\Dd_{Z/Z'}^{\bot }\left( \Omega _{Z/Z'}^{d+2},\omega_{Y/Y'}\right)} \arrow{s} \\
\node{\Ll_{Z/Z'}^{0}\tensor _{\Oo_{Z}}\Uu_{\Oo_{Z}}\left(\Ll_{Z/Z'}^{0}\right)} \arrow{s} \arrow{e} \node{\Dd_{Z/Z'}^{\bot}\left( \Omega _{Z/Z'}^{d+1},\omega _{Y/Y'}\right)} \arrow{s} \\
\node{\Uu_{\Oo_{Z}}\left( \Ll_{Z/Z'}^{0}\right)} \arrow{s} \arrow{e} \node{\Dd_{Z/Z'}^{\bot }\left( \Omega _{Z/Z'}^{d},\omega_{Y/Y'}\right)} \\
\node{\Oo_{Z}} 
\end{diagram}
\label{eq.312}
\end{equation}
(Here we change from left to right action of $\Uu_{\Oo_{Z}}\left(
\Ll_{Z/Z'}^{0}\right) $ in the standard way, that is, by extending
the $\left( -1\right) $-involution on $\Ll_{Z/Z'}^{0}$
multiplicatively.)

The horizontal maps in~(\ref{eq.312}) are given by Lie differentiation, followed by
contraction with elements of $\wedge ^{*}\Ll_{Z/Z'}^{0}$, followed by
multiplication on the left by the restriction class $\rho $. It is an easy
calculation in local coordinates (see~(\ref{eq.A4}) in the Appendix) to see that the
horizontal maps are surjective.  Thus, using the Poincar\'{e}-Birkoff-Witt
Theorem we have increasing filtrations
\begin{equation*}
\Dd_{Z/Z'}^{\bot }\left( \Omega _{Z/Z'}^{d+k},\omega
_{Y/Y'}\right) _{\left[ m\right] }\subseteq \Dd_{Z/Z'}^{\bot }\left( \Omega _{Z/Z'}^{d+k},\omega _{Y/Y'}\right)
\end{equation*}
and surjections
\begin{equation*}
\wedge ^{k}\Ll_{Z/Z'}^{0}\tensor _{\Oo_{Z}}S_{\Oo%
_{Z}}^{m}\left( \Ll_{Z/Z'}^{0}\right) \arrow{e} \frac{\Dd%
_{Z/Z'}^{\bot }\left( \Omega _{Z/Z'}^{d+k},\omega
_{Y/Y'}\right) _{\left[ m\right] }}{\Dd_{Z/Z'}^{\bot
}\left( \Omega _{Z/Z'}^{d+k},\omega _{Y/Y'}\right)
_{\left[ m-1\right] }}.
\end{equation*}
Using commutativity, skew commutativity, and the fact that
\begin{equation*}
\rho \cdot \left\< \xi \right. \left| {\ }\right\> =0
\end{equation*}
if $\xi $ points along the fibers of $p_{Y'}$, we see that
this last map factors through
\begin{equation}
\wedge ^{k}\Ll_{Z/Z'}^{0}\tensor _{\Oo_{Y}}S_{\Oo%
_{Z}}^{m}N_{Y\backslash Y'\times _{Z'}Z}\longrightarrow
\frac{\Dd_{Z/Z'}^{\bot }\left( \Omega _{Z/Z'}^{d+k},\omega _{Y/Y'}\right) _{\left[ m\right] }}{\Dd%
_{Z/Z'}^{\bot }\left( \Omega _{Z/Z'}^{d+k},\omega
_{Y/Y'}\right) _{\left[ m-1\right] }}
\label{eq.313}
\end{equation}
which must therefore be surjective.

Now the filtration $\left\{ \left( _{\left[ 0\right] }\Dd_{Y/Y'}^{Z/Z'}\right) _{\left[ m\right] }\right\} $ induced on $_{\left[
0\right] }\Dd_{Y/Y'}^{Z/Z'}$ from
\begin{equation*}
\Dd_{Z/Z'}^{\bot }\left( \Omega _{Z/Z'}^{d+1},\omega
_{Y/Y'}\right) _{\left[ m-1\right] }\circ \partial \subseteq \
_{\left[ 0\right] }\Dd_{Y/Y'}^{Z/Z'}
\end{equation*}
(with a degree shift of one) is the same one as that induced from the
standard filtration on the preimage of $\Dd_{Y/Y'}^{Z/Z'}$
in $\Uu_{\Oo_{Z}}\left( \Ll_{Z/Z'}^{0}\right)$
given by number of Lie differentiations. On the other hand,
let $\left\{ \left( \Uu_{Z/Z'}^{0}\right) _{\left[ m\right] }\right\}
$ denote the filtration on $\Uu_{Z/Z'}^{0}$ given by number of Lie
differentiations. Using~(\ref{eq.276}) and the commutative diagram
\begin{equation*}
\begin{diagram}
\node{\Uu_{Z/Z'}^{0}} \arrow{s,V} \arrow{e}  \node{_{\left[ 0\right] }\Dd_{Y/Y'}^{Z/Z'}=\ker \left( \circ \partial \right) } \arrow{s,V} \arrow{e}  \node{\Oo_{Y'}=\Bbb{C}} \arrow{s,V} \\
\node{\Uu_{\Oo_{Z}}\left( \Ll_{Z/Z'}^{0}\right) } \arrow{e}  \node{\Dd_{Z/Z'}^{\bot }\left( \Omega _{Z/Z'}^{d},\omega _{Y/Y'}\right) } \arrow{e}  \node{\Oo_{Y}} 
\end{diagram}
\end{equation*}
we see that $\left\{ \left( \Uu_{Z/Z'}^{0}\right) _{\left[ m\right]
}\right\} $ is subordinate to the filtration on $\Uu_{Z/Z'}^{0}$
induced from the standard one on $\Uu_{\Oo_{Z}}\left( \Ll_{Z/Z'}^{0}\right)$.
Also, if $m>0$, $\left( \Uu_{Z/Z'}^{0}\right)
_{\left[ m\right] }$ maps onto the $m$-th graded quotient of the latter
filtration since it surjects to $S_{\Oo_{Y}}^{m}N_{Y\backslash
Y'\times _{Z'}Z}$. Therefore by~(\ref{eq.313}) the map
\begin{eqnarray*}
\Uu_{Z/Z'}^{0} &\longrightarrow& \Dd_{Y/Y'}^{Z/Z'} \\
u &\longmapsto& \rho \cdot u 
\end{eqnarray*}
must also be surjective. Thus the only thing left to check is that the
natural surjection
\begin{equation*}
S_{\Oo_{Y}}^{m}N_{Y\backslash Y'\times _{Z'}Z}\arrow{e} \frac{\left( _{\left[ 0\right] }\Dd_{Y/Y'}^{Z/Z'}\right) _{\left[ m\right] }}{\left( _{\left[ 0\right] }%
\Dd_{Y/Y'}^{Z/Z'}\right) _{\left[ m-1\right] }}
\end{equation*}
is actually injective. But this is immediately seen by writing this map in
local coordinates. (See the Appendix.)\bigskip

\begin{lem}\label{lem.32}
\begin{enumerate}
\item[i)]{}\label{lem.32.i}
$\Dd_{Y/Y'}^{Z/Z'}$ is the $%
\Uu $-submodule of $\Dd_{Z/Z'}^{\bot }\left( \Omega
_{Z/Z'}^{d},\omega _{Y/Y'}\right) $ generated by $\rho $.
The natural filtration $\left\{ \Uu_{\left[ m\right] }\right\} $ on $\Uu$
induces an increasing filtration on $\Dd_{Y/Y'}^{Z/Z'} $ such that
\begin{equation*}
Gr_{0}\left( \Dd_{Y/Y'}^{Z/Z'}\right) =\Oo%
_{Y'}\cdot \rho ,
\end{equation*}
and, if $m>0$,
\begin{equation*}
Gr_{m}\left( \Dd_{Y/Y'}^{Z/Z'}\right)
\end{equation*}
is itself filtered with associated gradeds being given by
\begin{equation*}
S_{m,r}=S_{\Oo_{Y}}^{m-r}N_{Y\backslash Y'\times _{Z'}Z}\tensor _{\Oo_{Y'}}S_{\Oo_{Z'}}^{r}T_{Z'}.
\end{equation*}

\item[ii)]{}\label{lem.32.ii}
Let $Y_{0}$ denote a fiber of $p_{Y'}$. Suppose we
are given locally around a point $y\in Y_{0}$ a projection
\begin{equation*}
\begin{diagram}
\node{Z} \arrow{e,t}{\nu} \node{Y_{0}}
\end{diagram}
\end{equation*}
(where, abusing notation, we temporarily denote the neighborhood of $y\in Z$
again by the letter $Z$). Let $\left\{ \xi _{k}\right\} _{k=1,\ldots ,c+\dim
Z'}$ be a commuting set of vector fields in $\Ll_{Y',Z'}$ which locally frame the fibers of $\nu $. Then $\Dd%
_{Y/Y'}^{Z/Z'}$ is in fact locally generated over $\Oo%
_{Y'}\cdot \rho $ by the subalgebra of $\Uu$ generated by the vector
fields in the set $\left\{ \Oo_{Y_{0}}\cdot \xi _{k}\right\} $%
.
\end{enumerate}
\end{lem}

\textbf{Proof:} By~(\ref{eq.273}) we have an isomorphism
\begin{equation*}
\frac{_{\left[ r\right] }\Dd_{Y/Y'}^{Z/Z'}}{_{\left[
r-1\right] }\Dd_{Y/Y'}^{Z/Z'}}=\ _{\left[ 0\right] }%
\Dd_{Y/Y'}^{Z/Z'}\ \tensor _{\Oo_{Z'}}\;S_{\Oo_{Z'}}^{r}T_{Z'}.
\end{equation*}
So by Lemma~\ref{lem.31} we have inductively that $\Uu^{0}\cdot \Uu_{\left[ r\right]
}$ generates$_{\left[ r\right] }\Dd_{Y/Y'}^{Z/Z'} $.\bigskip

Besides $\rho _{Y}\cdot \Uu$ we have another useful subset of $\Dd\left(
\Omega _{Z/Z'},\Omega _{Y/Y'}\right) $ which
surjects to $\Dd_{Y/Y'}^{Z/Z'}$. One checks easily,
from the fact that $\rho $ is $\partial $-closed (see the Appendix) and the
basic expression
\begin{equation*}
L_{\varsigma }=\left\< \zeta \right. \left| \partial \left( {}\right)
\right\> +\partial \left\< \zeta \right. \left| {\ }\right\> ,
\end{equation*}
that, for $\zeta _{i}\in \Ll_{Y',Z}$,
\begin{equation*}
\rho \cdot L_{\varsigma _{1}}\cdot \ldots \cdot L_{\varsigma _{m}}=\rho
\cdot \left( \left\< \zeta _{1}\right. \left| {\ }\right\> \cdot
\partial \cdot \ldots \cdot \left\< \zeta _{m}\right. \left| {\ }\right\> \cdot \partial \right)
\end{equation*}
as elements of $\Dd_{Y/Y'}^{Z/Z'}$, where $%
\left\< \zeta \right. \left| {\ }\right\> $ is the contraction
operator with the vector field $\zeta $.

Another fact that will be important in what follows is:\bigskip
\begin{lem}\label{lem.34}
\begin{enumerate}
\item[i)]{}\label{lem.34.i}
Referring to \S\ref{sec.2.3}, let $\Jj'$ be a
sheaf of ideals on $X'$ annihilating the image of $Y'$.
Suppose that the natural map
\begin{equation*}
\Uu_{\Jj'}'\arrow{e} \sheafHom_{\Bbb{C}}\left(
\left( \Oo_{X'}/\Jj'\right) _{y'},%
\Bbb{C}\right)
\end{equation*}
is surjective for all $y'\in Y'$. Then for
\begin{equation*}
\left( \Dd_{Y/Y'}^{X/X'}\right) _{\Jj}=\mathcal{H%
}\left( \Dd\left( \Omega _{X_{\Jj}/X_{\Jj'}}^{\le
d},\Omega _{Y/Y'}\right) ^{*,*}\right) \subseteq \Dd%
_{Y/Y'}^{X/X'}
\end{equation*}
we have
\begin{equation*}
\rho \cdot \Uu_{\Jj}=\left( \Dd_{Y/Y'}^{X/X'}\right) _{\Jj}.
\end{equation*}

\item[ii)]{}\label{lem.34.ii}
Fix a foliation of a neighborhood of $y$ in $X$ so that leaves project
isomorphically to a neighborhood of $h'\left( p_{Y'}\left( y\right) \right) $ in $Z'$ and project to points on $%
Y'$. This foliation determines a unique lifting of $\Ll_{Y',Z'}$ into $\Ll_{Y',Z}$ and of enveloping
algebras
\begin{eqnarray*}
\Uu_{y'}' &\longrightarrow& \Uu_{y} \\
u' &\longmapsto& u=u\left( u'\right). 
\end{eqnarray*}
Let $\tilde{\Jj}'\subseteq \Jj'$ be a sub-ideal
such that
\begin{equation*}
\frac{\Jj'}{\tilde{\Jj}'}=\Bbb{C},
\end{equation*}
and let $x'\in \Jj'$ and $u'\in \Uu'$
be such that
\begin{equation*}
u'\left( x'\right) =1.
\end{equation*}
Then
\begin{equation*}
\rho _{Y}\cdot u\left( u'\right) \left( x'\cdot \eta
\right) =\rho _{Y}\left( \eta \right) ,\quad \eta \in \Omega
_{X/X'}^{p}.
\end{equation*}
\end{enumerate}
\end{lem}

\textbf{Proof:} The foliation gives a local product structure
\begin{equation*}
\begin{diagram}
\node{X} \arrow{e,t}{\left( p_{X'},\pi \right)} \node{X'\times Z_{0}.}
\end{diagram}
\end{equation*}
Furthermore $Z_{0}$ itself is locally a product
\begin{equation*}
\begin{diagram}
\node{Z_{0}} \arrow{e,t}{\left( \alpha ,\beta \right)} \node{Y_{0}\times S.}
\end{diagram}
\end{equation*}
Let $V_{S}$ be the locally defined subalgebra of $\Uu$ generated by the
lifting of a maximal commuting set of vector fields on $S$ via the local
product structure (see Lemma~\ref{lem.32}(ii)). $V_{S}$ then corresponds to a
uniquely chosen set of local coordinates on $S$ such that, if we write
differential forms in terms of these coordinates on $S$ and arbitraily
chosen coordinates on $Y_{0}$, the action of $V_{S}$ via Lie
differentiation is just given by the action of $V_{S}$ on the
coefficient functions. We use these coordinates as part of a local
coordinate system on
\begin{equation*}
X=X'\times Y_{0}\times S
\end{equation*}
given by the product structure and compute $\Dd_{Y/Y'}^{Z/Z'}$ locally as in the Appendix. By Lemma~\ref{lem.32}(ii)
\begin{equation*}
\beta ^{-1}V_{S}+u\left( \Uu'\right) _{y}
\end{equation*}
generates $\Dd_{Y/Y'}^{Z/Z'}$. An easy computation in
local coordinates shows that the action of $u\left( \Uu'\right) _{y}$
on
\begin{equation*}
\pi ^{*}\Omega _{Z_{0}}^{p}=\pi ^{*}\left( \sum\nolimits_{s}\alpha
^{*}\Omega _{Y_{0}}^{p-s}\tensor \beta ^{*}\Omega _{S}^{s}\right)
\end{equation*}
is simply given by the action of $\Uu_{y'}'$ on coefficient
functions and this action commutes with the action of $\left( \beta \circ
\pi \right) ^{-1}V_{S}$ . Both i) and ii) of the Lemma then follow
easily.\bigskip

{}From~(\ref{eq.275}) we also have
\begin{equation*}
\mathcal{H}\left( \Dd^{\bot }\left( \Omega _{Z/Z'}^{>d},\omega
_{Y/Y'}\right) ,\circ \partial _{d}\right) =\Dd^{\bot
}\left( \Omega _{Z/Z'}^{d+1},\omega _{Y/Y'}\right)
\circ \partial _{d}\subseteq \Dd_{Y/Y'}^{Z/Z'},
\end{equation*}
and from~\S\ref{sec.2.5} and Lemmas~(\ref{lem.31}) and~(\ref{lem.32}) we have:\bigskip

\begin{cor}\label{cor.35}
The ``Koszul complex''
\begin{eqnarray*}
	&\ldots& \arrow{e} \rho \cdot \left( \wedge ^{r+1}\Ll_{Y',Z}^{0}\right) \tensor
		\Uu\arrow{e,t}{\circ \partial }\rho \cdot
		\left( \wedge ^{r}\Ll_{Y',Z}^{0}\right) \tensor \Uu\arrow{e,t}{\circ \partial }\ldots \\
	&\ldots& \arrow{e,t}{\circ \partial }\rho
		\cdot \Ll_{Y',Z}^{0}\tensor \Uu\arrow{e,t}{\circ \partial }%
		\Dd^{\bot }\left( \Omega _{Z/Z'}^{d+1},\omega _{Y/Y'}\right) \circ \partial _{d} 
\end{eqnarray*}
is a resolution of $\Dd^{\bot }\left( \Omega _{Z/Z'}^{d+1},\omega _{Y/Y'}\right) \circ \partial _{d}$ as a
(right) $\Uu$-module and the map
\begin{equation*}
\rho \cdot \left( \wedge ^{\ge -1}\Ll_{Y',Z}^{0}\right) \tensor
\Uu\arrow{e} \Dd^{\bot }\left( \Omega _{Z/Z'}^{\ge d},\omega
_{Y/Y'}\right)
\end{equation*}
is an isomorphism. (Of course this is not a $\Uu$-free resolution because of
the relations implicit in the factor `` $\rho \cdot $ '' on the
left-hand-side.)
\end{cor}

\textbf{Proof:} The only thing we have to check is that we have not changed
anything in the proof of Lemma~\ref{lem.31} by replacement of $\Uu_{\Oo
_{Z}}\left( \Ll_{Z/Z'}^{0}\right) $ by the slightly smaller $\Uu$.
But as above we still get surjectivity on each graded quotient and therefore
surjectivity.

\newpage

\section{Canonical morphisms to $\Dd_{Y/Y'}^{Z/Z'}$\ associated to deformations}\label{sec.functor}

As we stated at the outset, given the situation of~\S\ref{sec.2.1} we are going to
make a correspondence between certain homomorphisms into $\Dd%
_{Y/Y'}^{Z/Z'}$ and deformations of $Y/Y'$. In
one direction that correspondence is essentially contained in the next
lemma.\bigskip

\begin{lem}\label{lem.41}
Suppose we have the commutative diagram
\begin{equation*}
\begin{diagram}
\node{Y} \arrow{s,r}{ p_{Y'}} \arrow{e,t}{i} \node{X} \arrow{s,r}{ p_{X'}} \arrow{e,t}{j} \node{Z} \arrow{s,r}{ p_{Z'}} \\
\node{Y'} \arrow{e,t}{i'} \node{X'} \arrow{e,t}{j'} \node{Z'} 
\end{diagram}
\end{equation*}
as in~\S\ref{sec.2.1}. Further suppose that the left-hand rectangle is a fibered product.
Then, whenever $y'=p_{Y'}\left( y\right) $, the above
diagram induces a natural map
\begin{equation*}
\Uu_{y'}'\arrow{e} \left( \Dd_{Y/Y'}^{Z/Z'}\right) _{y}.
\end{equation*}
\end{lem}

\textbf{Proof:} As in Lemma~\ref{lem.34} fix a foliation of a neighborhood of $y$
in $X$ so that leaves project isomorphically to a neighborhood of $h'\left( p_{Y'}\left( y\right) \right) $ in $Z'$ and
project to points on $Y'$. This foliation determines a unique
lifting of $\Ll_{Y',Z'}$ into $\Ll_{Y',Z}$
and of enveloping algebras
\begin{eqnarray*}
\Uu_{y'}' &\longrightarrow& \Uu_{y} \\
u' &\longmapsto& u=u\left( u'\right) 
\end{eqnarray*}
and so induces a map:
\begin{equation*}
\Uu_{y'}'\arrow{e} \left( \Dd_{Y/Y'}^{Z/Z'}\right) _{y}
\end{equation*}
To see that the maps are independent of the foliation, notice that, since
the maps from fibers of $Y/Y'$ to fibers of $Z/Z'$ are
local isomorphisms, any operator in the right ideal generated by a Lie
differentiation by an element of $\Ll_{Y,Z}^{0}$ annihilates $\rho $.
Independence then comes from~(\ref{eq.222}) and the fact that any two liftings of
an element of $\Ll_{Y',Z'}$ differ by an element of $%
\Ll_{Y',Z}^{0}$.\bigskip

Thus we have:\bigskip

\begin{cor}\label{cor.42}
Suppose with respect to~(\ref{eq.211}) we have a
commutative diagram
\begin{equation*}
\begin{diagram}
\node{Y} \arrow{s,r}{ p_{Y'}} \arrow{e}  \node{W} \arrow{s,r}{ q_{X'}} \arrow{e}  \node{X} \arrow{s,r}{ p_{X'}} \\
\node{Y'} \arrow{e}  \node{X'} \arrow{e}  \node{X^{\prime }} 
\end{diagram}
\end{equation*}
for which the left-hand rectangle is a fibered product. Then there is a canonical
induced morphism
\begin{equation*}
\Uu_{y'}'\arrow{e} \left( \Dd_{Y/Y'}^{Z/Z'}\right) _{y}.
\end{equation*}
\end{cor}

\textbf{Proof:} We apply Lemma~\ref{lem.41} to the case in which we replace $Z$
with $W$. Then use~\S\ref{sec.2.8}.\bigskip

Later we will also need:\bigskip

\begin{lem}\label{lem.43}
\begin{enumerate}
\item[i)]{}\label{lem.43.i}
Let $\Jj'$ be a sheaf of ideals on $%
X'$ which annihilates the image of $Y'$. Suppose that the
natural map
\begin{equation*}
\left( \Uu_{\Jj'}'\right) _{y'}\arrow{e}
\Hom_{\Bbb{C}}\left( \left( \Oo_{X'}/\Jj'\right) _{y'},\mathbf{C}\right)
\end{equation*}
is surjective for all $y'\in Y'$. Suppose in addition that
we have a commutative diagram
\begin{equation*}
\begin{diagram}
\node{Y} \arrow{s,r}{ p_{Y'}} \arrow{e}  \node{W_{J}} \arrow{s,r}{ q_{X'}} \arrow{e}  \node{X} \arrow{s,r}{ p_{X'}} \\
\node{Y'} \arrow{e}  \node{X_{\Jj'}'} \arrow{e}  \node{X'} 
\end{diagram}
\end{equation*}
and further suppose that the left-hand rectangle is a fibered product. Then there is
a canonical induced morphism
\begin{equation*}
\left( \Uu_{\Jj'}'\right) _{y'}\arrow{e}
\left( \Dd_{Y/Y'}^{Z/Z'}\right) _{y}.
\end{equation*}

\item[ii)]{}\label{lem.43.ii}
Suppose that $\tilde{\Jj}'\subseteq \Jj'$ is
a sub-ideal such that
\begin{equation*}
\frac{\Jj'}{\tilde{\Jj}'}=\Bbb{C}.
\end{equation*}
Suppose that we have two extensions
\begin{equation*}
\begin{diagram}
\node{Y} \arrow{s,r}{ p_{Y'}} \arrow{e}  \node{W^{i}} \arrow{s,r}{ q_{X'}^{i}} \arrow{e}  \node{X} \arrow{s,r}{ p_{X'}} \\
\node{Y'} \arrow{e}  \node{X'} \arrow{e,=} \node{X'} 
\end{diagram}
\end{equation*}
for $i=0,1$, which agree over $X_{\Jj'}'$. Suppose we
have foliations of $X$ as above which are compatible with $W^{0}$ and $W^{1}$
respectively and which agree over $X_{\Jj'}'$. Let
\begin{equation*}
u^{0},u^{1}\in \Uu
\end{equation*}
be the liftings of
\begin{equation*}
u'\in \Uu_{\tilde{\Jj}'}'
\end{equation*}
given in Lemma~\ref{lem.34}(ii) by the respective foliations. Then there is a vector
field
\begin{equation*}
\chi \in \Ll_{Y',Z}^{0}
\end{equation*}
such that, over $X_{\Jj'}'$%
\begin{equation*}
\rho \cdot \left( u^{1}-u^{0}\right) =\rho \cdot L_{\chi }.
\end{equation*}
\end{enumerate}
\end{lem}

\textbf{Proof:} i) Apply Corollary~\ref{cor.42} locally near $y$ in $Y$ to any
family which extends
\begin{equation*}
q_{X'}^{-1}X_{\Jj'}\arrow{e} X_{\Jj%
'}
\end{equation*}
to a family over $X'$. The map given by Corollary~\ref{cor.42} restricts
to a map on $\left( \Uu_{\Jj'}'\right) _{y'}$.
By Lemma~\ref{lem.34} the map does not depend on the choice of extension.

ii) With respect to the family $W^{0}$, we write everything in the local
coordinates chosen as in the proofs of Lemma~\ref{lem.34}. Suppose that, for $S$ as
in the proof of Lemma~\ref{lem.34} with coordinates $t_{1},\ldots ,t_{c}$%
, the extension
\begin{equation*}
W^{0}\arrow{e,t}{q_{X}}X
\end{equation*}
is given by the conditions $t_{k}=0$ . Let $x'\in \tilde{\Jj}'$ be a generator of $\frac{\Jj'}{\tilde{\Jj}%
'}$. The extension
\begin{equation*}
W^{1}\arrow{e,t}{q_{X}}X'
\end{equation*}
is given over $X_{\tilde{\Jj}'}$ by the equations
\begin{equation*}
\tilde{t}_{k}=0,
\end{equation*}
where
\begin{equation*}
\tilde{t}_{k}=t_{k}+a_{k}\left( y\right) x'.
\end{equation*}
Then we have:
\begin{eqnarray*}
\left( u^{1}-u^{0}\right) \left( t_{k}\right) &=& u^{1}\left( t_{k}\right) \\
	&=& u^{1}\left( \tilde{t}_{k}-a_{k}\left( y\right) x'\right) \\
	&=& -u^{1}\left( a_{k}\left( y\right) x'\right) \\
	&=& -a_{k}\left( y\right) \cdot u'\left( x'\right) 
\end{eqnarray*}
Also $\left( u^{1}-u^{0}\right) $ annihilates pull-backs of functions from $%
X'\times Y_{0}$ and all products $t_{k}\cdot t_{k'}\cdot g$. Thus:
\begin{equation*}
\rho \cdot \left( u^{1}-u^{0}\right) =\rho \cdot
L_{-\sum\nolimits_{k}a_{k}\left( y\right) u'\left( x'\right) \frac{\partial }{\partial t_{k}}}=\rho \cdot L_{-u'\sum\nolimits_{k}a_{k}\left( y\right) x'\frac{\partial }{\partial
t_{k}}}
\end{equation*}

\newpage

\section{Dolbeault resolution of $\Dd_{Y/Y'}^{Z/Z'} $}\label{sec.resol}

Again in this section we assume that~(\ref{eq.211}) is a relative immersion. From~(\ref{eq.274})--(\ref{eq.275}) we have pairings
\begin{equation}
\begin{diagram}
\node{R\left( p_{Y'}\right) _{*}\Dd\left( \Omega _{Z/Z'}^{>d},\Omega _{Y/Y'}^{*}\right) \tensor R\left( p_{Z'}\right) _{*}\left( \Omega _{Z/Z'}^{d+1}
	\rightarrow \ldots \rightarrow \Omega _{Z/Z'}^{n-1}\rightarrow \omega _{Z/Z'}\right)} \arrow{s} \\
\node{R\left( p_{Y'}\right) _{*}\left( \Omega _{Y/Y'}^{*}\right)}
\end{diagram}
\label{eq.511}
\end{equation}
and
\begin{equation}
\begin{diagram}
\node{R\left( p_{Y'}\right) _{*}\left( \Dd_{Y/Y'}^{Z/Z'}\right) \tensor R\left( p_{Z'}\right) _{*}\left( \Oo_{Z}\rightarrow \Omega _{Z/Z'}^{1}\rightarrow \ldots \rightarrow \Omega _{Z/Z'}^{d}\right)} \arrow{s} \\
\node{R\left( p_{Y'}\right) _{*}\left( \Omega _{Y/Y'}^{*}\right)}
\end{diagram}
\label{eq.512}
\end{equation}
Referring to~(\ref{eq.274}), these pairings are related by the commutative diagram:
\begin{equation}
\thinD
\begin{diagram}
\node{R\left( p_{Y'}\right) _{*}\left( \Dd_{Y/Y'}^{Z/Z'}\right)
		\tensor R\left( p_{Z'}\right)_{*}\left( \Omega _{Z/Z'}^{\le d}\right)}
	\arrow{se} \\
\node{R\left( p_{Y'}\right) _{*}\Dd\left( \Omega _{Z/Z'}^{>d}\left[ 1\right] ,\Omega _{Y/Y'}^{*}\right)
		\tensor R\left( p_{Z'}\right)_{*}\left( \Omega _{Z/Z'}^{\le d}\right)}
	\arrow{n,t}{(\circ \partial)\tensor 1}
	\arrow{s,b}{1\tensor\partial}
	\node{R\left(p_{Y'}\right) _{*}\left( \Omega _{Y/Y'}^{*}\right)} \\
\node{R\left( p_{Y'}\right) _{*}\Dd\left( \Omega _{Z/Z'}^{>d}\left[ 1\right] ,\Omega _{Y/Y'}^{*}\right)
		\tensor R\left(p_{Z'}\right) _{*}\left( \Omega _{Z/Z'}^{>d}\right) \left[ 1\right]}
	\arrow{ne} 
\end{diagram}
\label{eq.513}
\end{equation}
If, for example, $\Dd_{Y/Y'}^{Z/Z'}$ were a subsheaf
of $\Dd\left( \Omega _{Z/Z'}^{\le d},\Omega _{Y/Y'}^{*}\right) $, we could compute its pairings at the level of complexes of
sheaves, but in general we are not in this situation and the pairing is only
defined in the derived category. To work at the level of complexes of
sheaves it is useful to reformulate $\left( p_{Y'}\right)
_{*}\left( \Dd_{Y/Y'}^{Z/Z'}\right) $ in terms of
Dolbeault complexes.

Recall from~(\ref{eq.261}) that
\begin{equation*}
\Dd_{Y/Y'}^{Z/Z'}=\mathcal{H}\left(
\sum\nolimits_{p\ge c}D\left( \Omega _{Z/Z'}^{n-p},\Omega
_{Y/Y'}^{q}\right) ,\left( \circ \partial \right) +\left( -1\right)
^{p}\left( \partial \circ \right) \right) .
\end{equation*}
We let $\Aa_{Z/Z'}^{p,r}$, resp. $\Aa_{Y/Y'}^{q,s}$denote the sheaves of $C^{\infty }$-differential forms of type
$(p,r)$ on $Z/Z'$,resp. of type $(q,s)$ on $Y/Y'$. Now let $%
\Ff$ be a locally free $\Oo_{Z}$-module, and let $\Dd%
_{h}^{\infty }\left( \Ff,\Aa_{Y/Y'}^{q,s}\right) $ denote
the sheaf of $\Bbb{C}$-linear operators $T$ from the sheaf of (holomorphic)
sections of $\Ff$ to the sheaf of $C^{\infty }$-sections $\Aa_{Y/Y'}^{q,s}$
such that, for some $n$ and for arbitrary
holomorphic functions $f_{0},\ldots ,f_{n}$,
\begin{equation*}
\left[ \ldots \left[ \left[ T,f_{0}\right] ,f_{1}\right] \ldots
,f_{n}\right] =0.
\end{equation*}
We can consider $\Dd\left( \Ff,\Omega _{Y/Y'}^{*}\right) $
as a subsheaf of $\Dd_{h}^{\infty }\left( \Ff,\Aa_{Y/Y'}^{*,0}\right) $,
and conclude that, we have, over $Y'$:
\begin{equation}
\begin{array}{rcl}
\Oo_{Y}\tensor _{\Dd\left( \omega _{Y/Y'}\right) }\Dd\left( \Ff,\omega _{Y/Y'}\right)
	&\equiv& \mathcal{H}\left( \Dd\left( \Ff,\Omega _{Y/Y^{\prime }}^{*}\right) ,\left( \partial \circ \right) \right) \\
	&\equiv& \mathcal{H}\left( \Dd_{h}^{\infty }\left( \Ff,\Aa_{Y/Y'}^{*}\right) ,\left( d\circ \right) \right) 
\end{array}
\label{eq.521}
\end{equation}
The first equivalence is derived as in \S\S\ref{sec.2.4}--\ref{sec.2.5} or~\cite[13.2]{bib.B} and the
second is derived from the quasi-isomorphism
\begin{equation*}
\begin{diagram}
\node{\left( \Omega _{Y/Y'}^{*},\partial \right)} \arrow{e} \node{\left( \Aa_{Y/Y'}^{*},d\right) .}
\end{diagram}
\end{equation*}
(There is a small point here that the cohomology of the left-hand complex is
the sheaf of holomorphic sections of a holomorphic vector bundle whereas the
cohomology of the right-hand complex is the sheaf of $C^{\infty }$-sections
of a $C^{\infty }$-vector bundle. Our assertion is that the underlying
complex vector bundles are naturally isomorphic.)

Alternatively let
\begin{equation*}
\Dd^{\infty }\left( \Ff,\Aa_{Y/Y'}^{q,s}\right)
\end{equation*}
denote the sheaf of $C^{\infty }$-differential operators from the sheaf of $%
C^{\infty }$-sections of $\Ff$ to the sheaf of $C^{\infty }$-sections
of $\Aa_{Y/Y'}^{q,s}$. We can characterize $\Dd_{h}^{\infty
}\left( \Ff,\Aa_{Y/Y'}^{q,s}\right) $ as the subsheaf of $%
\Dd^{\infty }\left( \Ff,\Aa_{Y/Y'}^{q,s}\right) $
consisting of those operators $P$ such that, for every holomorphic section $%
\zeta $ of $\Ff$ and holomorphic function $f$ on $Z$,
\begin{equation*}
P\left( \bar{f}\cdot \zeta \right) =\bar{f}\cdot P\left( \zeta \right) .
\end{equation*}
Thus $\Dd_{h}^{\infty }\left( \Ff,\Aa_{Y/Y'}^{q,s}\right) $ is a right $\overline{\Oo}_{Z}$-module. Tensoring
over $\overline{\Oo}_{Z}$ with the exact sequence of locally free left $\overline{\Oo}_{Z}$-modules
\begin{equation*}
	0
	\arrow{e} \overline{\Dd_{Z}\left(\Omega _{Z/Z'}^{n},\Oo_{Z}\right) }
	\arrow{e,t}{\overline{\circ \partial }} \cdots
	\arrow{e,t}{\circ \partial } \overline{\Dd_{Z}\left(\Omega _{Z/Z'}^{1},\Oo_{Z}\right) }
	\arrow{e,t}{\overline{\circ \partial }} \overline{\Dd_{Z}\left( \Oo_{Z},\Oo_{Z}\right) }
	\arrow{e} \overline{\Oo_{Z}}
	\arrow{e} 0
\end{equation*}
we obtain
\begin{equation}
\mathcal{H}\left( \Dd^{\infty }\left( \Ff,\Aa_{Y/Y'}^{*,*}\right) ,\left( d\circ \right) \right) \equiv \mathcal{H}\left(
\sum\nolimits_{r,q,s}^{\infty }\Dd^{\infty }\left( \Aa%
_{Z/Z'}^{0,n-r}\left( \Ff\right) ,\Aa_{Y/Y'}^{q,s}\right) ,\left( \circ \bar{\partial }\right) +\left( -1\right)
^{r}\left( d\circ \right) \right)
\label{eq.522}
\end{equation}
We apply~(\ref{eq.521}) and~(\ref{eq.522}) to the case in which $\Ff$ is a graded
quotient of $\left( \sum\nolimits_{p\ge c}\Omega _{Z/Z'}^{n-p},\partial \right) $. We obtain:
\begin{equation}
\Dd^{\bot }\left( \Omega _{Z/Z'}^{n-p},\omega _{Y/Y'}\right) =\mathcal{H}\left( \sum\nolimits_{r,q,s}\Dd^{\infty
}\left( \Aa_{Z/Z'}^{n-p,n-r},\Aa_{Y/Y'}^{q,s}\right) ,\left( \circ \bar{\partial }\right) +\left( -1\right)
^{p+r}\left( d\circ \right) \right)
\label{eq.53}
\end{equation}
The standard spectral sequence argument then allows us to conclude that the
natural inclusion of complexes gives a quasi-isomorphism:
\begin{equation}
\Dd_{Y/Y'}^{Z/Z'}=\mathcal{H}\left(
\sum\nolimits_{p\ge c}\sum\nolimits_{q,r,s}\Dd^{\infty }\left( \Aa%
_{Z/Z'}^{n-p,n-r},\Aa_{Y/Y'}^{q,s}\right) ,\left(
\circ d\right) +\left( -1\right) ^{p+r}\left( d\circ \right) \right) .
\label{eq.54}
\end{equation}
By the acylicity of the sheaves in~(\ref{eq.53}) we have:
\begin{equation*}
R\left( p_{Y'}\right) _{*}\left( \Dd^{\bot }\left( \Omega
_{Z/Z'}^{n-p},\omega _{Y/Y'}\right) \right) =\mathcal{%
H}\left( \sum\nolimits_{r,q,s}^{\infty }\Dd^{\infty }\left(
\Aa_{Z/Z'}^{n-p,n-r},\Aa_{Y/Y'}^{q,s}\right) ,\left( \circ
\bar{\partial }\right) +\left( -1\right) ^{p+r}\left( d\circ \right) \right)
\end{equation*}
where
\begin{equation*}
\Aa^{*,*}=\left( p_{Y'}\right) _{*}\Aa^{*,*}
\end{equation*}
and by the acyclicity of the sheaves in and~(\ref{eq.54}) we have
\begin{equation*}
R\left( p_{Y'}\right) _{*}\left( \Dd_{Y/Y'}^{Z/Z'}\right) =\mathcal{H}\left( \sum\nolimits_{p\ge
c}\sum\nolimits_{q,r,s}\Dd^{\infty }\left( \Aa_{Z/Z'}^{n-p,n-r},\Aa_{Y/Y'}^{q,s}\right) ,\left( \circ d\right) +\left(
-1\right) ^{p+r}\left( d\circ \right) \right) .
\end{equation*}
Thus for any fixed $p$ and $q$ the natural pairing of $\bar{\partial }$-complexes
\begin{equation}
\sum\nolimits_{r,s}\Dd^{\infty }\left( \Aa_{Z/Z'}^{p,r},\Aa_{Y/Y'}^{q,s}\right) \tensor
\sum\nolimits_{r}\Aa_{Z/Z'}^{p,r}\arrow{e}
\sum\nolimits_{s}\Aa_{Y/Y'}^{q,s}
\label{eq.57}
\end{equation}
obtained by ``letting differential operators operate'' induces the pairings~(\ref{eq.511}) and~(\ref{eq.512})
on Dolbeault cohomology. Thus, using Corollary~\ref{cor.35},
the pairings~(\ref{eq.511}) and~(\ref{eq.512}) are induced from~(\ref{eq.57}) and the map
\begin{equation*}
\begin{diagram}
\node{\left( p_{Y'}\right) _{*}\left( \Aa_{Z/Z'}^{0,*}\tensor \wedge ^{j}\left( \Ll_{Y',Z}\right) \tensor \Uu\right)}
\arrow{e} \node{\Dd^{\infty }\left( \Aa_{Z/Z'}^{*+j,*},\Aa_{Y/Y'}^{*,*}\right) .}
\end{diagram}
\end{equation*}

\newpage

\section{The Gauss-Manin connection}\label{sec.gauss}

If we are in the situation of Corollary~\ref{cor.42}, that is. we have a geometric
extension
\begin{equation}
\begin{diagram}
\node{Y} \arrow{s,r}{ p_{Y'}} \arrow{e}  \node{W} \arrow{s,r}{ q_{X'}} \arrow{e}  \node{Z} \arrow{s,r}{ p_{Z'}} \\
\node{Y'} \arrow{e,t}{i'} \node{X'} \arrow{e,t}{j'} \node{Z'} 
\end{diagram}
\label{eq.61}
\end{equation}
of $Y/Y$', then~(\ref{eq.512}) and Corollary~\ref{cor.42} imply that we have an induced
pairing
\begin{equation*}
\Uu'\tensor j^{\prime *}R\left( p_{Z'}\right)
_{*}\left( \Oo_{Z}\arrow{e} \Omega _{Z/Z'}^{1}\arrow{e} \ldots \arrow{e} \Omega _{Z/Z'}^{d}\right)
\arrow{e} R\left( p_{Y'}\right) _{*}\left( \Omega
_{Y/Y'}^{*}\right) .
\end{equation*}
We next wish to see that this pairing is nothing more than the pull-back map
\begin{eqnarray*}
\lefteqn{R\left( p_{Z'}\right) _{*}\left( \Oo_{Z}\rightarrow
	\Omega _{Z/Z'}^{1}\rightarrow \ldots \rightarrow \Omega _{Z/Z'}^{d}\right)} \\
	&\rightarrow&
	R\left( q_{X'}\right)_{*}\left( \Oo_{W}\rightarrow \Omega _{W/X'}^{1}\rightarrow
		\ldots \rightarrow \Omega _{W/X'}^{d}\right) \\
	&=& R\left( q_{X'}\right) _{*}\left( \Bbb{C}_{W}\right) 
\end{eqnarray*}
followed by the Gauss-Manin connection on the family $W/X'$.\bigskip

\begin{dfn}\label{dfn.63}
Let $\xi '\in \Ll_{Y',Z'}$. By a $C^{\infty }$-lifting of type $\left( 1,0\right) $ of $\xi
'$ to $X/X'$ we mean a $C^{\infty }$-vector field $\xi $
on $X$ such that

\begin{enumerate}
\item[i)]
$\left[ \xi ,\bar{\xi}\right] =0$,

\item[ii)]
$\xi \left( f'\circ p_{X'}\right) =\xi '\left( f'\right) $ for every $C^{\infty }$-function $f'$
on $X'$.
\end{enumerate}
\end{dfn}

Given~(\ref{eq.61}) one constructs a $C^{\infty }$-projection of a tubular
neighborhood of $Y/Y'$ in $W/X'$ compatible with the
projection $X'\rightarrow Y'$ and such that the fibers are
complex holomorphic polydisks isomorphic to $Z'$. Any such
projection induces a unique $C^{\infty }$-lifting of each $\xi '\in
\Ll_{Y',Z'}$ to a vector field tangent to the fibers of
the projection and so a homomorphism of sheaves of Lie algebras
\begin{equation*}
\Ll_{Y',Z'}\arrow{e} \Aa_{X}^{0}\tensor
\Ll_{Y',Z}.
\end{equation*}
and therefore an induced mapping
\begin{equation}
\Uu'\arrow{e} \Dd^{\infty }\left( \Aa_{Z/Z'}^{\le
d,*},\Aa_{Y/Y'}^{*,*}\right) .
\label{eq.642}
\end{equation}
As in~\S\ref{sec.2.2} Lie differentiation commutes with the action of $\circ d$ and so
that $\Uu'$ acts on the complex $\Dd^{\infty }\left(
\Aa_{Z/Z'}^{\le d,*},\Aa_{Y/Y'}^{*,*}\right) $ and the
mapping~(\ref{eq.642}) is simply the image of that action on the distinguished
generator $\rho $. Under the standard Cech-to-Dolbeault isomorphism,~(\ref{eq.642})
is identified with the mapping given by Corollary~\ref{cor.42}. Thus we get that
the pairing~(\ref{eq.512}) is given by:
\begin{equation*}
\begin{diagram}
\node{\Uu'\tensor \Dd^{\infty }\left( \Aa_{Z/Z'}^{\le d,*},\Aa_{Y/Y'}^{*,*}\right)}
	\arrow{e} \node{\Uu'\tensor \Dd^{\infty }\left( \Aa_{W/X'}^{*,*},\Aa_{Y/Y'}^{*,*}\right)}
	\arrow{e,t}{\text{Gauss-Manin}} \node{\Aa_{Y/Y'}^{*,*}}
\end{diagram}
\end{equation*}

\newpage

\section{The deformation criterion}\label{sec.defcrit}

Suppose now that we have commutative diagram
\begin{equation}
\begin{diagram}
\node{Y} \arrow{s,r}{ p_{Y'}} \arrow{e,t}{h} \node{Z} \arrow{s,r}{ p_{Z'}} \\
\node{Y'} \arrow{e,t}{h'} \node{Z'} 
\end{diagram}
\label{eq.71}
\end{equation}
as in~(\ref{eq.211}) where $Y/Y'\rightarrow Z/Z'$ is relatively
immersive. On one hand we wish to consider extensions of~(\ref{eq.71}) to diagrams
of the form
\begin{equation}
\begin{diagram}
\node{Y} \arrow{s,r}{ p_{Y'}} \arrow{e}  \node{W} \arrow{s,r}{ q_{X'}} \arrow{e}  \node{X} \arrow{s,r}{ p_{X'}} \\
\node{Y'} \arrow{e,t}{i'} \node{X'} \arrow{e,=} \node{X'} 
\end{diagram}
\label{eq.72}
\end{equation}
for which the left-hand rectangle is a fibered product. For example, if $Y'$
is a one-point space,~(\ref{eq.72}) amounts to giving a deformation of the
submanifold $Y$ in $Z/Z'$.

We first apply Corollary~\ref{cor.42} to associate to~(\ref{eq.72}) a family of natural
morphisms
\begin{equation*}
\begin{diagram}
\node{\Uu_{y'}'} \arrow{e,t}{\varphi _{y',y}} \node{\left( \Dd_{Y/Y'}^{Z/Z'}\right) _{y}}
\end{diagram}
\end{equation*}
whenever $y\in Y_{y'}=p_{Y'}^{-1}\left( y'\right)
$, and so to define a morphism
\begin{equation}
\Uu'\arrow{e} \left( p_{Y'}\right) _{*}\left( \Dd%
_{Y/Y'}^{Z/Z'}\right) .
\label{eq.73}
\end{equation}

In the rest of this section, we shall be concerned with the opposite
direction. Namely, given a diagram~(\ref{eq.71}) and a morphism~(\ref{eq.73}):

\begin{enumerate}
\item
What are the conditions on the morphism such that it is induced by a
diagram~(\ref{eq.72})?

\item
Does the morphism~(\ref{eq.73}) completely determine the diagram~(\ref{eq.72}) that
induces it?
\end{enumerate}

To our state the result precisely, we need two definitions.\bigskip

\begin{dfn}\label{dfn.74}
If $y'=p_{Y'}\left(
y\right) $, an $\Oo_{X/X'}$-linear homomorphism of Lie
algebras of tangent vector fields
\begin{equation*}
\psi :\left( \Ll_{Y',Z'}\right) _{y'}\arrow{e} \left( \Ll_{Y',Z}\right) _{y}
\end{equation*}
is called a lifting homomorphism if, for all $\xi '\in \left(
\Ll_{Y',Z'}\right) _{y'}$:

\begin{enumerate}
\item[i)]
$\psi \left( g\cdot \xi '\right) =(g\circ p_{X'})\cdot \psi \left( \xi '\right) $ for all $g\in \Oo%
_{X'}$,

\item[ii)]
$\psi \left( \xi '\right) (g\circ p_{X'})=\xi
'(g)$ for all $g\in \Oo_{X'}$.
\end{enumerate}
\end{dfn}

\begin{dfn}\label{dfn.75}
Given $y\in Y,\;y'\in Y'$ with $%
y'=p_{Y'}\left( y\right) $, a mapping
\begin{equation*}
\begin{diagram}
\node{\Uu_{y'}'} \arrow{e,t}{\varphi _{y',y}} \node{\left( \Dd_{Y/Y'}^{Z/Z'}\right) _{y}}
\end{diagram}
\end{equation*}
is called almost-multiplicative if there exists a lifting homomorphism
\begin{equation*}
\psi :\left( \Ll_{Y',Z'}\right) _{y'}\arrow{e} \left( \Ll_{Y',Z}\right) _{y}
\end{equation*}
such that, for the induced homomorphism
\begin{equation*}
\begin{diagram}
\node{\Uu_{y'}'} \arrow{e,t}{\psi} \node{\Uu_{y}}
\end{diagram}
\end{equation*}
and for all $u'\in \Uu'$,
\begin{equation*}
\varphi _{y'}\left( u'\right) =\rho \cdot \psi \left(
u'\right) ,
\end{equation*}
where $\rho $ is the distinguished generator of $\Dd_{Y/Y'}^{Z/Z'}$. (See~(\ref{eq.276}).)
\end{dfn}

We call a map
\begin{equation*}
\begin{diagram}
\node{\Uu_{y'}'} \arrow{e,t}{\varphi _{y'}} \node{\left( \left(p_{Y'}\right) _{*}\left( \Dd_{Y/Y'}^{Z/Z'}\right) \right) _{y'}}
\end{diagram}
\end{equation*}
\textit{almost multiplicative} if the induced maps
\begin{equation*}
\begin{diagram}
\node{\Uu_{y'}'} \arrow{e,t}{\varphi _{y',y}} \node{\left( \Dd_{Y/Y'}^{Z/Z'}\right) _{y}}
\end{diagram}
\end{equation*}
are almost multiplicative whenever $y\in Y,\;y'\in Y'$
with $y'=p_{Y'}\left( y\right) $. In short, almost
multiplicative homomorphisms are those induced locally by
``cotransversals,'' that is, by local foliations $X$ transverse to the
fibers of $Y/Y'$. (The important point here is that the existence
of $\psi $ is required only locally; a global lifting of $\Ll_{Y',Z'}$ to $\left( p_{Y'}\right) _{*}\left(
\Ll_{Y',Z}\right) $ is usually not possible.)

Notice that the mapping~(\ref{eq.73}) derived from~(\ref{eq.72}) and Corollary~\ref{cor.42} is, by
construction, almost multiplicative. By Lemma~\ref{lem.32}(ii), these almost
multiplicative maps are determined by

\begin{enumerate}
\item[i)]
choosing a commuting set of elements of $\Ll_{Y',Z'}$
which locally frame the tangent space of $Z'$ along $Y'$,

\item[ii)]
lifting these to a commuting family of vector fields on $X$.
\end{enumerate}

\noindent The chief result we wish to prove in \S\ref{sec.defcrit} can now be stated as
follows.\bigskip

\begin{thm}\label{thm.main}
There exists a natural one-to-one correspondence
between formal extensions of the diagram~(\ref{eq.71}) to diagrams~(\ref{eq.72}) and almost
multiplicative homomorphisms
\begin{equation*}
\begin{diagram}
\node{\Uu'} \arrow{e,t}{\varphi} \node{\left( p_{Y'}\right) _{*}\left( \Dd_{Y/Y'}^{Z/Z'}\right) .}
\end{diagram}
\end{equation*}
Here ``formal'' means that $W$ is a formal neighborhood of $i\left( Y\right)
$ in $X$.
\end{thm}

\textbf{Proof:} As already mentioned, the mapping
\begin{equation*}
\begin{diagram}
\node{\Uu'} \arrow{e,t}{\varphi} \node{\left( p_{Y'}\right)_{*}\left( \Dd_{Y/Y'}^{Z/Z'}\right)}
\end{diagram}
\end{equation*}
associated to a diagram~(\ref{eq.72}) is almost multiplicative by construction.
Also, since all morphisms are left $\Oo_{Y'}$-linear, it
suffices to establish the result in the case in which $Y'=\left\{
0\right\} $ is a point, which we will therefore assume throughout the proof.
So to complete the proof of the theorem it will suffice to establish the
following assertion:
The mapping
\begin{equation}
\begin{diagram}
\node{\Uu_{0}'} \arrow{e,t}{\varphi _{y}} \node{\left( \Dd_{Y/Y'}^{Z/Z'}\right) _{y}}
\end{diagram}
\label{eq.761}
\end{equation}
uniquely determines a (formal) local family
\begin{equation*}
W\arrow{e} Z'
\end{equation*}
fitting into~(\ref{eq.72}).

As is standard in local deformation theory, we work inductively over
increasing large (infinitesimal) neighborhoods of $\left\{ 0\right\} $ in $%
Z'$. We can extend the notion of almost multiplicative morphism to
the restriction of $\varphi _{y}$ to a scheme containing $\left\{ 0\right\} $
and lying in $Z'$ as follows:

Let $\Jj_{1}'$ be the ideal defining $\left\{ 0\right\} $ in $%
Z'$, and suppose we have a filtration
\begin{equation*}
\Jj_{s+1}'\subseteq \ldots \subseteq \Jj_{1}'\subseteq \Jj_{0}'=\Oo_{Z'}
\end{equation*}
by ideals on $Z'$ such that successive quotients
\begin{equation*}
J_{r}=\frac{\Jj_{r}'}{\Jj_{r+1}'}
\end{equation*}
are dimension-one complex vector spaces generated by
\begin{equation*}
x_{r}'\in \Jj_{r}'.
\end{equation*}
Then the set $\left\{ x_{j}'\right\} _{j=0,\ldots ,r}$ gives each
quotient $\Oo_{r}=\frac{\Oo_{Z'}}{J_{r+1}}$ the structure
of a $\Bbb{C}$-vector space. Let
\begin{equation*}
\begin{diagram}
\node{\Oo_{r}} \arrow{e,t}{D_{r}} \node{J_{r}=\Bbb{C}}
\end{diagram}
\end{equation*}
be the natural projection homomorphism associated with this basis. It is
immediate to see that $D_{r}$ is a differential operator of order $\le r$
with values in $\Bbb{C}$. Thus $D_{r}$ comes from an element
\begin{equation*}
u_{r}'\in \Uu_{0}'.
\end{equation*}
In fact, if we let
\begin{equation*}
\Uu_{\left[ r\right] }'=\left\{ u'\in \Uu':u'\left( g\right) \left( 0\right) =0\ for\ all\ g\in \Jj_{r+1}'\right\} ,
\end{equation*}
we have
\begin{equation*}
u_{r}'\in \Uu_{\left[ r\right] }'\subseteq \Uu'.
\end{equation*}

Letting $Z_{\left[ r\right] }'$ denote the scheme defined in $%
Z'$ by $\Jj_{r+1}'$, and referring to~(\ref{eq.233}), our
definition for ``infinitesimal almost-multiplicativity over $Z_{\left[
r\right] }'$'' is obtained from Definition~\ref{dfn.75} by replacing the
condition
\begin{equation*}
\varphi _{y}\left( u'\right) =\rho \cdot \psi \left( u'\right)
\end{equation*}
for $u\in \Uu'$ with the same condition restricted to $u'\in
\Uu_{\left[ r\right] }'$.

If we have constructed
\begin{equation*}
W_{\left[ s\right] }\arrow{e,t}{q}Z_{\left[ s\right]
}'
\end{equation*}
over $Z_{\left[ s\right] }'$, then by Lemma~\ref{lem.43}(i) we have a
uniquely determined almost multiplicative map
\begin{equation*}
\begin{diagram}
\node{\Uu_{\left[ s\right] }'} \arrow{e,t}{\varphi _{y}^{s}} \node{\left( \Dd_{Y/Y'}^{Z/Z'}\right) _{y}.}
\end{diagram}
\end{equation*}
We must show a natural bijective correspondence between extensions of $%
\varphi _{y}^{s}$ to an almost multiplicative map
\begin{equation*}
\begin{diagram}
\node{\Uu_{\left[ s+1\right] }'} \arrow{e,t}{\varphi _{y}^{s+1}} \node{\left( \Dd_{Y/Y'}^{Z/Z'}\right) _{y}}
\end{diagram}
\end{equation*}
and extensions of $W_{\left[ s\right] }\arrow{e,t}{q}%
Z_{\left[ s\right] }'$ to a family over $Z_{\left[ s+1\right]
}'$.

We let $Y_{0}$ denote the fiber of $p_{Y'}$. Let
\begin{equation*}
\Jj=\text{ideal sheaf of $Y_{0}$ in $Z$}
\end{equation*}
whose restriction gives the ideal sheaf $\Jj_{0}$ of $Y_{0}$ in $Z_{0}$
. At $y\in Y_{0}$ let
\begin{equation*}
\left\{ z_{k}\in \Jj_{y}\right\} _{k}
\end{equation*}
restrict to a free set of generators of $\Jj_{0}$. For $y_{0}\in Y_{0}$%
, abuse notation slightly and let $Y_{0}$ denote the local germ of $Y_{0}$
at $y_{0}$. Standard deformation theory shows that the set of local
extensions of a (flat) family
\begin{equation*}
\begin{diagram}
\node{Y_{0}} \arrow{s,r}{ p_{Y'}} \arrow{e}  \node{W_{\left[ s\right] }} \arrow{s,r}{q} \arrow{e}  \node{Z} \arrow{s,r}{ p_{Z'}} \\
\node{\left\{ 0\right\}} \arrow{e,V} \node{Z_{\left[ s\right] }'} \arrow{e} \node{Z'} 
\end{diagram}
\end{equation*}
over $Z_{\left[ s\right] }'$ to a (flat) family over $Z_{\left[
s+1\right] }'$ is a principal homogeneous space for
\begin{equation*}
\Hom\left( \frac{\Jj_{0}}{\Jj_{0}^{2}},p_{Y'}^{*}\frac{%
\Jj_{s+1}'}{\Jj_{s+2}'}\right) _{y}=\Hom\left(
\frac{\Jj_{0}}{\Jj_{0}^{2}},p_{Y'}^{*}J_{s+1}\right) _{y}
\end{equation*}
(See~\cite[Chapter I]{bib.Ko}.) More precisely, fix a local product structure
\begin{equation*}
Z=Y_{0}\times S\times Z'
\end{equation*}
such that $p_{Z'}$ is given by projection onto the last
factor. Suppose that functions
\begin{equation*}
z_{k}=w_{k}+a_{k,0}\left( y\right) x_{0}'+\ldots +a_{k,s}\left(
y\right) x_{s}'
\end{equation*}
are a free set of local generators of the ideal of $W_{\left[ s\right] }$ in
\begin{equation*}
Z_{\left[ s\right] }=Z_{\left[ s\right] }'\times _{Z'}Z.
\end{equation*}
(We use the same notation for functions on $Z'$ and their pullbacks
to $Z$.) Fix reference extensions of $z_{k}$ to $Z_{\left[ s+1\right] }$,
and continue to call them $z_{k}$. The possible choices of extensions of $%
W_{\left[ s\right] }$ to $W_{\left[ s+1\right] }$ are then uniquely
determined by the functions
\begin{equation}
\tilde{z}_{k}=z_{k}+a_{k}\left( y\right) x_{s+1}'.
\label{eq.762}
\end{equation}
Thus extensions are in one-to-one correspondence with elements of
\begin{equation*}
\Hom\left( \frac{\Jj_{0}}{\Jj_{0}^{2}},p_{Y'}^{*}\Jj%
_{s+1}'\right) _{y}=\left( N_{Y_{0}\backslash Z_{0}}\right) _{y}.
\end{equation*}
So we need only show that the almost multiplicative extension
\begin{equation*}
\begin{diagram}
\node{\Uu_{\left[ s+1\right] }'} \arrow{e,t}{\varphi _{y}^{s+1}} \node{\left( \Dd_{Y/Y'}^{Z/Z'}\right) _{y}}
\end{diagram}
\end{equation*}
induced via Lemma~\ref{lem.43}(i) by a chosen extension $W_{\left[ s+1\right] }$
determines the element of
\begin{equation*}
\left( N_{Y_{0}\backslash Z_{0}}\right) _{y}
\end{equation*}
that defines $W_{\left[ s+1\right] }$ with respect to the reference
extension. But by Lemma~\ref{lem.43}(ii) we conclude that
\begin{equation}
\rho \cdot \tilde{\psi}\left( u_{s+1}'\right) -\rho \cdot \psi
\left( u_{s+1}'\right) =\rho \cdot L_{-\sum\nolimits_{k}a_{k}\left(
y\right) \frac{\partial }{\partial w_{k}}}.
\label{eq.763}
\end{equation}

Thus we read off from~(\ref{eq.763}) the exact element of
\begin{equation*}
\left( N_{Y_{0}\backslash Z_{0}}\right) _{y}
\end{equation*}
that determined the extension of the deformation given by the equations $%
\tilde{z}_{k}=0$ in~(\ref{eq.762}). So the extension
\begin{equation*}
\begin{diagram}
\node{\Uu_{\left[ s+1\right] }'} \arrow{e,t}{\varphi _{y}^{s+1}} \node{\left( \Dd_{Y/Y'}^{Z/Z'}\right) _{y}}
\end{diagram}
\end{equation*}
determines the extension
\begin{equation*}
W_{\left[ s+1\right] }\arrow{e,t}{q}Z_{\left[ s+1\right]
}'.
\end{equation*}

But Lemma~\ref{lem.43}(i) tell us that the extension
\begin{equation*}
W_{\left[ s+1\right] }\arrow{e,t}{q}Z_{\left[ s+1\right]
}'
\end{equation*}
determines the extension
\begin{equation*}
\begin{diagram} \node{\Uu_{\left[ s+1\right] }'} \arrow{e,t}{\varphi _{y}^{s+1}} \node{ \left( \Dd_{Y/Y'}^{Z/Z'}\right) _{y}} \end{diagram}
\end{equation*}
Thus the proof of assertion~(\ref{eq.761}) is complete, which in turn completes the
proof of Theorem~\ref{thm.main}.\bigskip

Inside complex projective varieties, formal subvarieties give rise to actual
geometric subvarieties, so we have:\bigskip

\begin{thm}\label{thm.77}
Suppose that we have a diagram~(\ref{eq.71}) with the same
hypotheses as given there. Suppose, in addition, $p_{Y'}$ is
proper and $Z/Z'$ is relatively quasi-projective. Given an almost
multiplicative morphism
\begin{equation*}
\begin{diagram} \node{\Uu'} \arrow{e,t}{\varphi } \node{ \left( p_{Y'}\right) _{*}\left( \Dd_{Y/Y'}^{Z/Z'}\right) ,} \end{diagram}
\end{equation*}
there is a geometric deformation
\begin{equation*}
\begin{diagram}
\node{Y\;} \arrow{s,r}{ p_{Y'}} \arrow{e,V} \node{W} \arrow{s,r}{ q_{X'}} \arrow{e,V} \node{X} \arrow{s,r}{ p_{X'}} \\
\node{Y'} \arrow{e,t}{i'} \node{X'} \arrow{e,=} \node{X'} 
\end{diagram}
\end{equation*}
such that the associated
\begin{equation*}
\begin{diagram} \node{\Uu'} \arrow{e,t}{\tilde{\varphi}} \node{ \left( p_{Y'}\right) _{*}\left( \Dd_{Y/Y'}^{Z/Z'}\right)} \end{diagram}
\end{equation*}
approximates the morphism $\varphi $ to arbitrary pregiven order.
\end{thm}

Theorem~\ref{thm.main} gives a formal family
\begin{equation*}
\begin{diagram}
\node{Y} \arrow{s,r}{ p_{Y'}} \arrow{e,V} \node{\hat{W}} \arrow{s,r}{ q_{X'}} \arrow{e,V} \node{X} \arrow{s,r}{ p_{X'}} \\
\node{Y'} \arrow{e,t}{i} \node{\hat{X}'} \arrow{e,=} \node{X'} 
\end{diagram}
\end{equation*}
projective over $\hat{X}$. So the formal subvariety $\hat{W}$ can be
approximated to arbitrary order by a $W$ which is projective
over $X$~\cite[Theorem 1.4]{bib.A}, so over a sufficiently small analytic neighborhood of $Y$ in $%
X $.\bigskip

{}From the proof of Theorem~\ref{thm.main} and Lemma~\ref{lem.43}(ii) we have:\bigskip

\begin{cor}\label{cor.78}
\begin{enumerate}
\item[i)]{}\label{cor.78.i}
Let
\begin{equation*}
\tilde{\Jj}'\subseteq \Jj'
\end{equation*}
be sheaves of ideals on $X'$ contained in the ideal of $Y'$
such that
\begin{equation*}
\frac{\Jj'}{\tilde{\Jj}'}=\Oo_{Y'}.
\end{equation*}
If $\sheafHom$ denotes almost multiplicative homomorphisms of sheaves and $%
\varphi $ is an almost multiplicative homomorphism giving the family
\begin{equation*}
W_{\Jj'}\arrow{e,t}{q_{\Jj'}}X_{%
\Jj'},
\end{equation*}
the obstruction to extending the family to
\begin{equation*}
W_{\tilde{\Jj}'}\arrow{e,t}{q_{\tilde{\Jj}'}}X_{\tilde{\Jj}'}
\end{equation*}
is the element of
\begin{equation*}
\sheafHom_{\Oo_{Y'}}\left( \left( \frac{\Jj'}{\tilde{\Jj}'}\right) ^{*},R^{1}\left( p_{Y'}\right) _{*}N_{Y\backslash Y'\times _{Z'}Z}\right) =R^{1}\left( p_{Y'}\right)
_{*}N_{Y\backslash Y'\times _{Z'}Z}
\end{equation*}
which measures the obstruction to lifting $\varphi $ under the homomorphism
\begin{equation*}
\sheafHom\left( \Uu_{\tilde{\Jj}'}',\left(
p_{Y'}\right) _{*}\left( \Dd_{Y/Y'}^{Z/Z'}\right) \right) \arrow{e} \sheafHom\left( \Uu_{\Jj'}',\left( p_{Y'}\right) _{*}\left( \Dd
_{Y/Y'}^{Z/Z'}\right) \right) .
\end{equation*}

\item[ii)]{}\label{cor.78.ii}
The map
\begin{eqnarray*}
N_{Y\backslash Y'\times _{Z'}Z}
	&\longrightarrow& \Dd^{\bot }\left( \Omega _{Z/Z'}^{d+1},\omega _{Y/Y^{\prime }}\right) \\
\xi &\longmapsto& \left\< \xi \right. \left| {\ }\right\> 
\end{eqnarray*}
induces a natural morphism
\begin{equation*}
R^{1}\left( p_{Y'}\right) _{*}N_{Y\backslash Y'\times
_{Z'}Z}\arrow{e} R^{1}\left( p_{Y'}\right)
_{*}\left( \Dd^{\bot }\left( \Omega _{Z/Z'}^{d+1},\omega
_{Y/Y'}\right) \circ \partial _{d}\right)
\end{equation*}
so that, from i) we obtain an obstruction element
\begin{equation*}
\upsilon \in R^{1}\left( p_{Y'}\right) _{*}\left( \Dd%
^{\bot }\left( \Omega _{Z/Z'}^{d+1},\omega _{Y/Y'}\right) \circ \partial _{d}\right) .
\end{equation*}
\end{enumerate}
\end{cor}

\textbf{Proof:}
i) By Lemma~\ref{lem.43}(ii), the sheaf sequence
$0\rightarrow \sheafHom_{\Oo_{Y'}}\left( p_{Y'}^{*}\left( \frac{\Jj'}{\tilde{\Jj}'}\right)
^{*},N_{Y\backslash Y'\times _{Z'}Z}\right)
\rightarrow \sheafHom\left( \Uu_{\tilde{\Jj}'}',
\Dd_{Y/Y'}^{Z/Z'}\right) \rightarrow \sheafHom\left( \Uu_{%
\Jj'}',\Dd_{Y/Y'}^{Z/Z'}\right) \rightarrow 0$ is exact. The obstruction element is computed from a
covering $\left\{ V_{a}\right\} $ of $Z$ by coordinate disks as the Cech
cochain
\begin{equation*}
\chi _{ab}=-\sum\nolimits_{k}a_{k}^{ab}\left( y\right) \frac{\partial }{%
\partial v_{k}^{ab}}\in \Gamma _{\Uu_{a}\cap \Uu_{b}}\Ll_{Y',Z}^{0}
\end{equation*}
constructed as in the proof of Theorem~\ref{thm.main}.

ii) The second assertion then follows from the exact sequence~(\ref{eq.275}).

\newpage

\section{Some immediate consequences}\label{sec.cor}

\begin{cor}\label{cor.81}
A compact submanifold $Y$ of dimension $d$ has a
geometric deformation in a projective manifold $Z$ of dimension $n$ if and
only if there is an almost multiplicative homomorphism
\begin{equation*}
\begin{diagram} \node{\Dd'} \arrow{e,t}{\varphi } \node{ H^{0}\left( Y;\Dd_{Y/\left\{ 0\right\} }^{Z/\left\{ 0\right\} }\right)} \end{diagram}
\end{equation*}
where $\Dd'$ is the ring of constant-coefficient differential
operators in one variable and $\Delta $ is the unit disc. Thus, referring to~(\ref{eq.252}),
a necessary condition for $Y$ to have a deformation in $Z$ is that
\begin{equation*}
\underset{\leftarrow s}{\lim }H^{d}\left( \Omega _{Z}^{d}/\Ii_{Y,Z}^{s}\Omega _{Z}^{d}\right)
\end{equation*}
be infinite-dimensional.
\end{cor}

\textbf{Proof:} The inverse limit $\underset{\leftarrow }{\lim }
H^{d}\left( \Omega _{Z}^{d}/\Ii_{Y,Z}^{s}\Omega _{Z}^{d}\right) $
is the dual of
\begin{equation*}
H_{Y}^{c}\left( \Omega _{Z}^{c}\right) \cong H^{0}\left( \Dd^{\bot
}\left( \Omega _{Z}^{d},\omega _{Z}\right) \right) .
\end{equation*}
(See~\cite[Theorem 7.13]{bib.B}.)\bigskip

Again suppose $Y=\left\{ 0\right\}$.  For a geometric family~(\ref{eq.72}), it is
easy to see, by induction on the order of neighborhood of $0\in Z'$%
, that
\begin{equation*}
\begin{diagram} \node{\Uu_{\left\{ 0\right\} }'} \arrow{e,t}{\varphi } \node{ H^{0}\left( D_{Y_{0}/\left\{ 0\right\} }^{W/X'}\right)} \end{diagram}
\end{equation*}
is in fact an isomorphism since that is trivially true for the associated
graded objects with respect to the filtration in Lemma~\ref{lem.32}. It is
natural to ask whether the composition
\begin{equation*}
\begin{diagram} \node{\Dd_{\left\{ 0\right\} }^{X'}} \arrow{e,t}{\varphi } \node{ H^{0}\left( \Dd_{Y_{0}/\left\{ 0\right\} }^{Z/Z'}\right)} \end{diagram}
\end{equation*}
associated to a diagram~(\ref{eq.72}) is an isomorphism if the associated
first-order map
\begin{equation*}
\begin{diagram} \node{T_{X',0}} \arrow{e,t}{\varphi } \node{ H^{0}\left( N_{Y_{0}\backslash Z}\right)} \end{diagram}
\end{equation*}
is. This is not in general true. The simplest example of this is perhaps the
classical logarithmic transform
\begin{equation*}
Z=Z_{0}=\frac{\Bbb{C}\times E}{\left\{ \left( t,e\right) =\left(
-t,e+h\right) \right\} },\quad Z'=\left\{ 0\right\}
\end{equation*}
where $E$ is an elliptic curve, $h$ is a non-trivial half-period, and $%
Y_{0}$ is given by $t=0$. Then
\begin{equation*}
H^{0}\left( N_{Y_{0}\backslash Z_{0}}\right) =H^{1}\left( N_{Y_{0}\backslash
Z_{0}}\right) =0
\end{equation*}
but the local isomorphisms
\begin{equation*}
\begin{diagram}
\node{\left( \Dd_{\left\{ 0\right\} \times E/\left\{ 0\right\} }^{\Bbb{C}\times E/\Bbb{C}}\right) _{\left( 0,e\right) }}
	\arrow{e,t}{\approx} \node{\left( \Dd_{\left\{ 0\right\} \times E/\left\{ 0\right\} }^{\Bbb{C}\times E/\Bbb{C} }\right) _{\left( 0,e+h\right) }}
	\arrow{e,t}{\approx} \node{\left( \Dd_{Y/\left\{ 0\right\} }^{X_{0}/\left\{ 0\right\} }\right) _{\left( 0,e\right) }} \\
\node{\rho \cdot \frac{\partial }{\partial t}}
	\arrow{e,T} \node{\rho \cdot \left( -\frac{\partial }{\partial t}\right)}
	\arrow{e,T} \node{\rho \cdot \xi} 
\end{diagram}
\end{equation*}
imply that there is a non-trivial second-order element of
\begin{equation*}
H^{0}\left( \Dd_{Y_{0}/\left\{ 0\right\} }^{Z/\left\{ 0\right\} }\right)
\end{equation*}
given everywhere locally by
\begin{equation*}
\rho \cdot \xi \cdot \xi .
\end{equation*}

However, there are situations in which, for example, the map
\begin{equation*}
\varphi :\left( \Dd_{\left\{ 0\right\} /\left\{ 0\right\} }^{X'/\left\{ 0\right\} }\right) \arrow{e} \left( p_{Y'}\right)
_{*}\left( \Dd_{\left\{ 0\right\} /\left\{ 0\right\} }^{X'\times X_{0}/\left\{ 0\right\} }\right) =H^{0}\left( \Dd_{\left\{
0\right\} /\left\{ 0\right\} }^{X_{0}/\left\{ 0\right\} }\right)
\end{equation*}
induced by a geometric family~(\ref{eq.72}) is in fact an isomorphism.\bigskip

\begin{lem}\label{lem.82}
For a geometric family~(\ref{eq.72}), suppose that the natural
maps
\begin{equation*}
\begin{diagram} \node{S^{m}T_{X',0}} \arrow{e,t}{\varphi _{m}} \node{ H^{0}\left( S^{m}N_{Y_{0}\backslash Z_{0}}\right)} \end{diagram}
\end{equation*}
are isomorphisms for each $m$, that is
\begin{equation*}
T_{X',0}\arrow{e,t}{\cong }H^{0}\left(
N_{Y_{0}\backslash Z_{0}}\right)
\end{equation*}
and
\begin{equation*}
S^{m}H^{0}\left( N_{Y_{0}\backslash Z_{0}}\right) =H^{0}\left(
S^{m}N_{Y_{0}\backslash Z_{0}}\right) .
\end{equation*}
(This is the case if, for example, $N_{Y_{0}\backslash Z_{0}}=\left( \Oo%
_{Y_{0}}\right) ^{c'}$.) Then the natural map
\begin{equation*}
\begin{diagram} \node{\Dd_{\left\{ 0\right\} }^{X'}} \arrow{e,t}{\varphi } \node{ H^{0}\left( \Dd_{Y_{0}/\left\{ 0\right\} }^{Z_{0}/\left\{ 0\right\} }\right)} \end{diagram}
\end{equation*}
is an isomorphism.
\end{lem}

\textbf{Proof:} As in the proof of Theorem~\ref{thm.main}, we proceed by induction on
the order of the neighborhood of $0$ in $Z'$. Using Lemma~\ref{lem.31},
the lemma is then an immediate consequence of the commutative diagram with
exact rows:
\begin{equation*}
\setlength{\dgARROWLENGTH}{1em}
\begin{diagram}
\node{0} \arrow{e} \node{\left( \Dd_{\left\{ 0\right\} }^{X'}\right)^{s}} \arrow{s,=}
	\arrow[2]{e} \node[2]{\left( \Dd_{\left\{ 0\right\} }^{X'}\right) ^{s+1}} \arrow{s}
	\arrow[2]{e} \node[2]{S^{s+1}T_{X',0}} \arrow{s,=} \arrow{e} \node{0} \\
\node{0} \arrow{e} \node{H^{0}\left( \left( \Dd_{Y_{0}/\left\{ 0\right\}}^{Z_{0}/\left\{ 0\right\} }\right)^{s}\right)}
	\arrow[2]{e} \node[2]{H^{0}\left(\left( \Dd_{Y_{0}/\left\{ 0\right\} }^{Z_{0}/\left\{ 0\right\} }\right)^{s+1}\right)}
	\arrow[2]{e} \node[2]{S^{s+1}N_{Y_{0}\backslash Z_{0}}} 
\end{diagram}
\end{equation*}
\bigskip

We now turn to some more general consequences of Theorem~\ref{thm.main}. First,
suppose that in~(\ref{eq.72}) we have a proper smooth family:
\begin{equation*}
\begin{diagram}
\node{W} \arrow{s,r}{ q_{X'}} \arrow{e}  \node{Z} \arrow{s,r}{ p_{X'}} \\
\node{X'} \arrow{e,t}{j'} \node{Z'} 
\end{diagram}
\end{equation*}
We would like to have a criterion to decide when the fibers $W_{x'}$ of $q_{X'}$ are of constant modulus (that is,
mutually isomorphic) when $x'$ moves along $\left\{ y'\right\} \times Z'$. Referring to~\S\ref{sec.2.6} we let
\begin{equation*}
\Zz_{Y/Y'}^{Z/Z'}=\ker \left( \Dd\left( \Omega
_{Z/Z'}^{\le d},\Omega _{Y/Y'}^{*}\right)
\overset{\left(\circ \partial \right) \pm \left( \partial \circ \right)}{\longrightarrow}
\Dd\left( \Omega _{Z/Z'}^{\le d},\Omega _{Y/Y'}^{*}\right)
\right) .
\end{equation*}
\bigskip

\begin{thm}\label{thm.83}
The family
\begin{equation*}
\begin{diagram}
\node{Y} \arrow{s,r}{ p_{Y'}} \arrow{e,V} \node{W} \arrow{s,r}{ q_{X'}} \arrow{e,V} \node{X'\times _{Z'}Z} \arrow{s,r}{ p_{X'}} \\
\node{Y'} \arrow{e,t}{i'} \node{X'} \arrow{e,=} \node{X'} 
\end{diagram}
\end{equation*}
in~(\ref{eq.72}) is expressible as a family
\begin{equation*}
\begin{diagram}
\node{Y} \arrow{s,r}{ p_{Y'}} \arrow{e,V} \node{X'\times_{Z'}Y} \arrow{s,r}{ q_{X'}} \arrow{e,t}{} \node{X'\times_{Z'}Z} \arrow{s,r}{ p_{X'}} \\
\node{Y'} \arrow{e,t}{i'} \node{X'} \arrow{e,=} \node{X'} 
\end{diagram}
\end{equation*}
if and only if the associated almost multiplicative morphism
\begin{equation*}
\begin{diagram} \node{\Uu'} \arrow{e,t}{\varphi } \node{ \left( p_{Y'}\right) _{*}\left( \Dd_{Y/Y'}^{Z/Z'}\right)} \end{diagram}
\end{equation*}
lifts to an almost multiplicative map
\begin{equation*}
\begin{diagram} \node{\Uu'} \arrow{e,t}{\tilde{\varphi}} \node{ \left( p_{Y'}\right) _{*}\left( \Zz_{Y/Y'}^{Z/Z'}\right) .} \end{diagram}
\end{equation*}
(Here the definition of almost multiplicative is extended in the obvious way
for a map to the right $\Uu$-module $\Zz_{Y/Y'}^{Z/Z'}$%
.)
\end{thm}

\textbf{Proof:} Suppose first that we have a lifting
\begin{equation*}
\begin{diagram} \node{\Uu'} \arrow{e,t}{\tilde{\varphi}} \node{ \left( p_{Y'}\right) _{*}\left( \Zz_{Y/Y'}^{Z/Z'}\right) .} \end{diagram}
\end{equation*}
Replace Z in~(\ref{eq.71}) and~(\ref{eq.72}) with $Y\times Z$ and the map
\begin{equation*}
Y\arrow{e,t}{h}Z
\end{equation*}
with:
\begin{eqnarray*}
Y &\overset{\tilde{h}}{\longrightarrow}& Y\times Z \\
y &\longmapsto& \left( y,h\left( y\right) \right) 
\end{eqnarray*}
The decomposition
\begin{equation*}
\Omega _{Y\times Z}^{\le d}=\sum\nolimits_{i+j\le d}p_{Y}^{*}\Omega
_{Y}^{i}\tensor p_{Z}^{*}\Omega _{Z}^{j}
\end{equation*}
induces a natural map
\begin{equation*}
\Dd\left( \Omega _{Z}^{\le d},\Omega _{Y}^{*}\right) \arrow{e}
\Dd\left( \Omega _{Y\times Z}^{\le d},\Omega _{Y}^{*}\right)
\end{equation*}
which in turn induces
\begin{equation*}
\left( p_{Y'}\right) _{*}\left( \Zz_{Y/Y'}^{Z/Z'}\right) \arrow{e} \left( p_{Y'}\right)
_{*}\left( \Dd_{Y/Y'}^{Y\times Z/Z'}\right) .
\end{equation*}
To see that the composition
\begin{equation*}
\hat{\varphi}:
\Uu'\arrow{e} \left( p_{Y'}\right) _{*}\left( \Zz_{Y/Y'}^{Z/Z'}\right) \arrow{e} \left( p_{Y'}\right) _{*}\left( \Dd_{Y/Y^{\prime }}^{Y\times Z/Z'}\right)
\end{equation*}
is almost multiplicative, notice that it is given by lifting vector fields
so that they always are subordinate to the foliation of $Y\times X$ with
leaves $\left\{ point\right\} \times X$ . Now apply Theorem~\ref{thm.main} to the
almost multiplicative morphism $\hat{\varphi}$. In the other direction, the
product structure allows a globally consistent choice of the local vector
fields $\psi \left( \xi '\right) $ in Definition~\ref{dfn.75}.

\newpage

\section{The obstruction pairing}\label{sec.obst}

Suppose again that we are given the relative immersion
\begin{equation*}
\begin{diagram}
\node{Y} \arrow{s,r}{ p_{Y'}} \arrow{e,t}{h} \node{Z\;} \arrow{s,r}{ p_{Z'}} \\
\node{Y'} \arrow{e,t}{h'} \node{Z'} 
\end{diagram}
\end{equation*}
and that $p_{Z'}$ is proper and relatively K\"{a}hler. Suppose
now that we are in the situation of Corollary~\ref{cor.78} where we have an
extension
\begin{equation*}
W_{\Jj'}\arrow{e,t}{q_{\Jj'}}%
X'
\end{equation*}
and $\tilde{\Jj}'\subseteq \Jj'$ with $\frac{%
\Jj'}{\tilde{\Jj}'}=\Oo_{Y'}.$

Suppose further that $\left( \frac{\Jj'}{\tilde{\Jj}%
'}\right) ^{\vee }$ is spanned by
\begin{equation}
u'=v'\cdot \xi '\in \Uu_{\tilde{\Jj}'}'
\label{eq.92}
\end{equation}
for some $v'\in \Uu_{\Jj'}'$ and some vector
field $\xi '\in \Ll_{Y',Z'}$. We can then
compute a relation that must be satisfied by the obstruction element $%
\upsilon $ given by Corollary~\ref{cor.78}(ii). This computation will be carried out
in this final section of this paper.

We apply the pairing~(\ref{eq.513}) to a Cech cochain
\begin{equation*}
\left\{ \chi _{ab}\in \left( p_{Y'}\right) _{*}\left(
\left. N_{Y\backslash Y'\times _{Z'}Z}\right|
_{V_{a}\cap V_{b}}\right) \right\}
\end{equation*}
giving the element $\upsilon $ in Corollary~\ref{cor.78}, where, as in Lemma~\ref{lem.34},
we have an open cover $\left\{ X_{a}\right\} $ of $X$ such that, for
each $X_{a}$ we have local decompositions
\begin{equation*}
X_{a}:=Y_{0}\times S\times Y'\times Z'
\end{equation*}
which locally simultaneously trivialize the families $Y/Y'$ and $%
Z/Z'$. On $X_{a}\cap X_{b}$, the product structure
associated to each of the two trivializations determines a lifting of
\begin{equation*}
\xi '\in \Ll_{Y',Z'}.
\end{equation*}
The difference of the two liftings gives
\begin{equation*}
\theta _{ab}\left( \xi '\right) \in \left( p_{X'}\right) _{*}\left( \left. \Ll_{Y',Z}^{0}\right| _{X_{a}\cap
X_{b}}\right) .
\end{equation*}
We denote by
\begin{equation*}
\begin{diagram} \node{\Uu'} \arrow{e,t}{u_{a}} \node{ \left. \Uu\right| _{X_{a}}} \end{diagram}
\end{equation*}
the (multiplicative) lifting morphism corresponding to the trivialization of
$X_{a}$. By Corollary~\ref{cor.78}
\begin{eqnarray*}
\rho \cdot L_{\chi _{ab}}
	&=& \rho \cdot \left( u_{b}\left( u^{\prime }\right) -u_{a}\left( u'\right) \right) \\
	&=& \rho \cdot \left( u_{b}\left( v'\cdot \xi '\right) -u_{a}\left( v'\cdot \xi '\right) \right) \\
	&=& \rho \cdot \left( u_{b}\left( v'\right) \cdot u_{b}\left( \xi '\right) -u_{a}\left( v^{\prime }\right) \cdot u_{b}\left( \xi '\right) \right) \\
	&=& \rho \cdot u\left( v'\right) \cdot \theta _{ab}\left( \xi '\right) 
\end{eqnarray*}
since
\begin{equation*}
\rho \cdot u\left( v'\right) =\rho \cdot u_{a}\left( v'\right) =\rho \cdot u_{b}\left( v'\right)
\end{equation*}
because liftings to $W_{\Jj'}$ are already well defined,
independently of local trivializations. Thus, from~(\ref{eq.513}) we have:\bigskip

\begin{thm}\label{thm.93}
Under the decomposition assumption~(\ref{eq.92}), the map
\begin{equation*}
R^{m}\left( p_{Z'}\right) _{*}\left( \Omega _{Z/Z'}^{>d}\right) \arrow{e,t}{\upsilon \cdot }R^{m+1}\left(
p_{Y'}\right) _{*}\left( \Omega _{Y/Y'}^{*}\right)
\end{equation*}
given by the pairing~(\ref{eq.511}) evaluated on the obstruction class $\upsilon $
is the same as the map
\begin{equation*}
\begin{diagram}
\node{R^{m}\left( p_{Z'}\right) _{*}\left( \Omega _{Z/Z^{\prime }}^{>d}\right)}
	\arrow{e} \node{R^{m+1}\left( p_{Z'}\right) _{*}\left( \Omega _{Z/Z'}^{d}\right)}
	\arrow{e} \node{R^{m+1}\left( p_{Y'}\right) _{*}\left( \Omega _{Y/Y'}^{*}\right) } \\
\node{\eta} \arrow{e,T} \node{\left\< \theta _{ab}\right. \left| {\ }\right\> \cup \eta}
	\arrow{e,T} \node{u\left( v'\right) \left( \left\< \theta _{ab}\right.  \left| {\ }\right\> \cup \eta \right) }
\end{diagram}
\end{equation*}
where the first map in the composition is the Kodaira-Spencer map for the
family $Z/Z'$ and the second is given by the pairing~(\ref{eq.512})
(induced from the Gauss-Manin connection--see \S\ref{sec.gauss}).\bigskip
\end{thm}

\begin{cor}\label{cor.94}
Assume the factorization~(\ref{eq.92}). Also assume that $%
p_{Y'}$ can be factored as

\begin{equation*}
Y\arrow{e,t}{\pi }V'\arrow{e,t}{\sigma }%
Y'
\end{equation*}
with $\pi $ proper, smooth, and of relative dimension $d'\le d$ so
that the extension $q_{\Jj'}$ can be factored as
\begin{equation*}
W_{\Jj'}'\arrow{e,t}{\pi }V_{\Jj%
'}'\arrow{e,t}{\sigma }X_{\Jj'}'.
\end{equation*}
Then a particular case of the obstruction map in Theorem~\ref{thm.93} becomes
\begin{equation*}
\begin{diagram}
\node{R^{d-1}\left( p_{Z'}\right) _{*}\left( \Omega _{Z/Z^{\prime }}^{d+1}\right) } \arrow{s,r}{\text{Kodaira-Spencer}} \\
\node{R^{d}\left( p_{Z'}\right) _{*}\left( \Omega _{Z/Z^{\prime }}^{d}\right) } \arrow{s,r}{\text{pull-back}} \\
\node{R^{d}\left( q_{\Jj'}\right) _{*}\left( \omega _{W_{\Jj'}'/X_{'\Jj'}}\right) } \arrow{s} \\
\node{\sigma _{*}\left( \omega _{V_{\Jj'}'/X_{\Jj'}'}\tensor R^{d'}\left( \pi \right)_{*}\left( \omega _{W_{\Jj'}'/V_{\Jj^{\prime }}'}\right) \right) }
	\arrow{s,r}{\text{Gauss-Manin}} \\
\node{\sigma _{*}\left( \omega _{V'/Y'}\tensor R^{d}\left( \pi \right) _{*}\left( \omega _{Y/V'}\right) \right) }
	\arrow{s,r}{\text{integration over fiber}} \\
\node{\sigma _{*}\omega _{V'/Y'}}
\end{diagram}
\end{equation*}
\end{cor}

We first apply Corollary~\ref{cor.94} in the case $d=d'$, that is, $%
p_{Y'}$ itself is proper.\bigskip

\def\acks{\cite{bib.R},~\cite{bib.Ka}}
\begin{cor}[\acks]\label{cor.95}
Suppose in Corollary~\ref{cor.94} we have $%
d=d'$ and are given that the cohomology class
\begin{equation*}
\gamma :=c(Y/Y')\in R^{2c}\left( p_{Z'}\right)
_{*}\left( \Bbb{C}_{Z}\right)
\end{equation*}
stays of type $\left( c,c\right) $. (Recall that $c=n-d$ is the codimension
of the fiber of $Y/Y'$ in the corresponding fiber of $Z/Z'$%
.) Then the obstruction pairing
\begin{equation*}
R^{d-1}\left( p_{Z'}\right) _{*}\left( \Omega _{Z/Z'}^{d+1}\right) \arrow{e} \Oo_{Y'}
\end{equation*}
vanishes.
\end{cor}

\textbf{Proof:} For $\omega \in R^{d-1}\left( p_{Z'}\right)
_{*}\left( \Omega _{Z/Z'}^{d+1}\right) $ and $\xi '\in
\Ll_{Y',Z'}$,
\begin{equation*}
\left. \left( \xi '\left( \omega \right) \right) \cdot \gamma
\right| _{W_{\Jj'}}=\left. \left\< \theta
_{ab}\right. \left| {\ }\right\> \cup \omega \right| _{W_{\Jj%
'}}\in R^{d}\left( q_{W_{\Jj'}}\right)
_{*}\left( \Omega _{W_{\Jj'}/X_{\Jj'}'}^{d}\right) =\Oo_{X_{\Jj'}'}
\end{equation*}
is zero because of the assumption that $\gamma $ is everywhere of type $%
\left( c,c\right) $ . Now apply Corollary~\ref{cor.94}.\bigskip

Finally we wish to generalize this application of Corollary~\ref{cor.94}. Suppose
as above that we have an extension
\begin{equation*}
\begin{diagram}
\node{Y} \arrow{s,r}{ p_{Y'}} \arrow{e}  \node{W_{\Jj'}} \arrow{s,r}{ q_{\Jj^{\prime }}} \arrow{e,V} \node{X} \arrow{s,r}{ p_{X'}} \\
\node{Y'} \arrow{e}  \node{X_{\Jj'}} \arrow{e,V} \node{X'} 
\end{diagram}
\end{equation*}
and wish to examine the obstruction pairing to extending over $X_{\tilde{\Jj}'}$.
(In particular we continue to assume the
factorization~(\ref{eq.92}).) In general, we have an Abel-Jacobi map in homology:
\begin{equation}
\begin{array}{rcl}
L_{d-d'}\sigma _{*}\left( \Bbb{Q}_{V_{\Jj^{\prime}}}\right)
		=\left( R^{d-d'}\sigma _{*}\left( \Bbb{Q}_{V_{\Jj'}}\right) \right) ^{\vee }
	&\longrightarrow& \left. \left( j^{\prime -1}\left( R^{d+d'}\left( p_{Z'}\right)_{*}
		\left( \Bbb{Q}_{Z}\right) \right) ^{\vee }\right) \right| _{X_{\Jj'}'} \\
\gamma
	&\longmapsto& \pi ^{-1}\left( \gamma \right) 
\end{array}
\label{eq.961}
\end{equation}
The image of~(\ref{eq.961}) must lie in (the dual of) a Hodge sub-structure
\begin{equation}
A_{\Jj'}\subseteq \left. \left( j^{\prime -1}\left(
R^{d+d'}\left( p_{Z'}\right) _{*}\left( \Bbb{Q}%
_{Z}\right) \right) \right) \right| _{X_{\Jj'}'}
\label{eq.962}
\end{equation}
of type
\begin{equation*}
\left( d,d'\right) +\ldots +\left( d',d\right) .
\end{equation*}
\bigskip

\begin{thm}\label{thm.97}
Assume that there is a rational sub-Hodge structure
of
\begin{equation*}
R^{d+d'}\left( p_{Z'}\right) _{*}\left( \Bbb{Q}%
_{Z}\right)
\end{equation*}
of type
\begin{equation*}
\left( d,d'\right) +\ldots +\left( d',d\right)
\end{equation*}
extending the one given by the image of the mapping~(\ref{eq.962}) (via the
Gauss-Manin connection). Then the image of the obstruction mapping in
Corollary~\ref{cor.94} lies in the subspace of $\sigma _{*}\left( \omega
_{V'/Y'}\right) $ given by the (relatively) exact
forms, that is, this subspace maps to zero under the natural map
\begin{equation*}
\sigma _{*}\left( \omega _{V'/Y'}\right)
\arrow{e} R^{d-d'}\sigma _{*}\left( \Bbb{C}_{V'}\right) .
\end{equation*}
\end{thm}

\textbf{Proof:} For $\omega \in R^{d-1}\left( p_{Z'}\right)
_{*}\left( \Omega _{Z/Z'}^{d+1}\right) $, we have by Corollary~\ref{cor.95} that
\begin{equation*}
\int\nolimits_{\pi ^{-1}\left( \gamma \right) }\left\< \upsilon \right.
\left| \omega \right\> =\int\nolimits_{\pi ^{-1}\left( \gamma \right)
}u\left( v'\right) \xi '\left( \omega \right) =v'\int\nolimits_{\pi ^{-1}\left( \gamma \right) }\xi '\left( \omega
\right) .
\end{equation*}
But the hypothesis that the sub-Hodge structure extends over$Z$ implies that
\begin{equation*}
\int\nolimits_{\pi ^{-1}\left( \gamma \right) }\xi '\left( \omega
\right) \equiv 0
\end{equation*}
for all $\gamma \in \left( R^{d-d'}\sigma _{*}\left( \Bbb{Q}%
_{V'}\right) \right) ^{\vee }$, $\omega \in
j^{-1}R^{d-1}\left( p_{Z'}\right) _{*}\left( \Omega
_{Z/Z'}^{d+1}\right) $ and $\xi '\in \Ll_{Y',Z'}$.\bigskip

Unfortunately, we do not, as yet, know any good applications of Theorem~\ref{thm.97}
to deformation theory. However there are situations, for example, the
case
\begin{equation*}
\begin{tabular}{c}
$n=3,\quad d=2,\quad d'=1$ \\
$\omega _{Z/Z'}=\Oo_{Z/Z'}$
\end{tabular}
\end{equation*}
in which the map in Corollary~\ref{cor.94} seems to carry most the obstruction theory
(see~\cite{bib.C} and~\cite{bib.CK}). So hopefully the above formalism will eventually yield
generalizations of variational Hodge-conjecture results like Theorem~4.13 of~\cite{bib.C}.

\appendix
\newpage

\section{Local coordinates}\label{sec.local}

Let $Y_{0}$ be a $d$-dimensional complex submanifold of the complex $n$-manifold
$Z_{0}$. For any $y\in Y_{0}$, choose a complete set of local
coordinates near $y$ in $Z_{0}$ made up as follows:

\begin{quote}
$x_{k}$ such that $x_{k}=0$ for all $k$ define $Y_{0}$ in $Z_{0}$,

$y_{j}$ such that $y_{j}=0$ for all $j$ define $y$ in $Y_{0}$.
\end{quote}

For any tuple of non-negative integers $K$ indexed by the set of indices $%
\{k\}$, we have associated the differential operator $\partial _{x}^{K}$
which give the composition of differentiation with respect each of the
coordinates $x_{k}$ the number of times given by the $k$-th integer in $K$.
Then the elements of the stalk of $\Dd_{Z_{0}/\left\{ 0\right\} }^{\bot
}\left( \Omega _{Z_{0}/\left\{ 0\right\} }^{p},\omega _{Y_{0}/\left\{
0\right\} }\right) $ can be uniquely written as
\begin{equation}
\sum\nolimits_{J,K,K',\left| J\right| +\left| K\right| =p}\left\{
a_{K,J,K'}\left( y\right) \left( \left( dx_{K}dy_{J}\right)
^{*}dy\right) \partial _{x}^{K'}\right\}
\label{eq.A1}
\end{equation}
where $J$ = subset of $\{1,\ldots ,d\}$, ordered by $<$, $K$ = subset of $%
\{1,\ldots ,c\}$, ordered by $<$,
\begin{equation*}
\begin{tabular}{l}
$dy_{J}=dx_{j_{1}}\wedge \ldots \wedge dx_{j_{\left| J\right| }}$ \\
$dx_{K}=dx_{k_{1}}\wedge \ldots \wedge dx_{k_{\left| K\right| }}$
\end{tabular}
\end{equation*}
and
\begin{equation*}
dy=dy_{1}\wedge \ldots \wedge y_{d}.
\end{equation*}

For example, the restriction map $\rho $ given above becomes $dy^{*}\cdot dy$%
. Notice next that right multiplication by $\partial _{x}^{K'}$ is
identical to right multiplication by the operator $L_{K'}$ given by
the composition of (commuting) Lie differentiations with respect to the
corresponding $|K'|$ vector fields $\frac{\partial }{\partial
x_{k'}}$.

Right multiplication of~(\ref{eq.A1}) by $\partial =\partial _{Z_{0}}$ is given by
composition with
\begin{equation*}
\sum \left( dx_{k}\wedge \right) \frac{\partial }{\partial x_{k}}+\sum
\left( dy_{j}\wedge \right) \frac{\partial }{\partial y_{j}},
\end{equation*}
and this action commutes with right multiplication by $\partial
_{x}^{K'}$, so that it suffices to study its action on elements
\begin{equation}
a_{J,K}(y)((dx_{K}\wedge dy_{J})^{*}\cdot dy).
\label{eq.A2}
\end{equation}

An easy computation gives that multiplying~(\ref{eq.A2}) on the right by $\partial $
gives
\begin{equation}
\begin{tabular}{l}
$\left\{ \left( -1\right) ^{\left| K\right| -1}\sum\limits_{j\in J}\frac{\sgn\left( j\left( J-\left\{ j\right\} \right) \right)
\partial a_{K,J}\left(y\right) }{\partial y_{j}}\left( \left( dx_{K}dy_{J-\left\{ j\right\}}\right) ^{*}dy\right) \right\} $ \\
$+\left\{ \sum\limits_{k\in K}\sgn\left( k\left( K-\left\{ k\right\} \right) \right)
a_{J,K}\left( y\right) \left( \left( dx_{K-\left\{ k\right\}}dy_{J}\right) ^{*}dy\right) \frac{\partial }{\partial x_{k}}\right\} $
\end{tabular}
\label{eq.A3}
\end{equation}
Thus, if $\partial $ is to annihilate a sum of elements~(\ref{eq.A3}), $a_{J,K}(y)=0$
whenever $K$ can be written as $K'\cup \{k\}$, that is, whenever $K$
is not empty. So we are left to compute $\partial (a(y)(dy^{*}\cdot dy))$,
which gives
\begin{equation*}
\sum\limits_{j=1,\ldots ,d}\frac{\pm \partial a\left( y\right) }{\partial
y_{j}}\left( \left( dy_{J-\left\{ j\right\} }\right) ^{*}dy\right) .
\end{equation*}
Thus, to get $\partial (a(y)(dy^{*}\cdot dy))=0$, $a(y)$ must be constant in
the $y$-coordinates.

To see that
\begin{equation}
\rho \cdot \wedge ^{r}\Ll_{Z_{0}}\cdot \Uu_{Z_{0}}=\Dd_{Z_{0}/\left\{
0\right\} }^{\bot }\left( \Omega _{Z_{0}/\left\{ 0\right\} }^{d+r},\omega
_{Y_{0}/\left\{ 0\right\} }\right)
\label{eq.A4}
\end{equation}
where $\Uu_{Z_{0}}$ is the enveloping algebra of the Lie algebra $\Ll_{Z_{0}}$,
we need only show how to generate elements
\begin{equation*}
\left( dx_{K}dy_{J}\right) ^{*}dy
\end{equation*}
by Lie differentiations and contractions. But this is immediate from the
fact that, for example,
\begin{equation*}
L_{y_{j}\frac{\partial }{\partial x_{k}}}
\end{equation*}
removes $dx_{k}$ and replaces it with $dy_{j}$.

\newpage
\bibliographystyle{amsalpha}
\bibliography{ndo}

\end{document}